# Inverse modeling of hydrologic systems with adaptive multi-fidelity Markov chain Monte Carlo simulations

Jiangjiang Zhang[1], Jun Man[1], Guang Lin[2], Laosheng Wu[3], and Lingzao Zeng[1]

[1] Zhejiang Provincial Key Laboratory of Agricultural Resources and Environment, Institute of Soil and Water Resources and Environmental Science, College of Environmental and Resource Sciences, Zhejiang University, Hangzhou, China,

[2] Department of Mathematics and School of Mechanical Engineering, Purdue University, West Lafayette, IN, USA,

[3] Department of Environmental Sciences, University of California, Riverside, CA, USA

Correspondence to:
L. Zeng,
lingzao@zju.edu.cn

# Abstract

Markov chain Monte Carlo (MCMC) simulation methods are widely used to assess parametric uncertainties of hydrologic models conditioned on measurements of observable state variables. However, when the model is CPU-intensive and high-dimensional, the computational cost of MCMC simulation will be prohibitive. In this situation, a CPU-efficient while less accurate low-fidelity model (e.g., a numerical model with a coarser discretization, or a data-driven surrogate) is usually adopted. Nowadays, multi-fidelity simulation methods that can take advantage of both the efficiency of the low-fidelity model and the accuracy of the high-fidelity model are gaining popularity. In the MCMC simulation, as the posterior distribution of the unknown model parameters is the region of interest, it is wise to distribute most of the computational budget (i.e., the high-fidelity model evaluations) therein. Based on this idea, in this paper we propose an adaptive multi-fidelity MCMC algorithm for efficient inverse modeling of hydrologic systems. In this method, we evaluate the high-fidelity model mainly in the posterior region through iteratively running MCMC based on a Gaussian process (GP) system that is adaptively constructed with multi-fidelity simulation. The error of the GP system is rigorously considered in the MCMC simulation and gradually reduced to a negligible level in the posterior region. Thus, the proposed method can obtain an accurate estimate of the posterior distribution with a small number of the high-fidelity model evaluations. The performance of the proposed method is demonstrated by three numerical case studies in inverse modeling of hydrologic systems.

# 1. Introduction

For a better understanding and management of hydrologic systems, there is a growing interest in applying numerical modeling techniques to conduct qualitative and quantitative analyses (*Anderson et al.*, 2015; *Vieux*, 2001). However, the existence of uncertainties in model structure, model parameters, initial and boundary conditions, measurement data, etc., would hinder the predictive accuracy of hydrologic modeling (*Clark et al.*, 2011; *Refsgaard et al.*, 2012; *Wagener & Gupta*, 2005). To reduce the predictive uncertainty of the hydrologic system of concern, it is common practice to calibrate the conceptual model against measurements of some state variables, e.g., hydraulic head, solute concentration, temperature and streamflow, through solving an inverse problem (*Hu et al.*, 2017; *Kang et al.*, 2017; *Zha et al.*, 2018; *Zhu et al.*, 2017).

When handling uncertainties in hydrologic modeling, inverse methods based on Bayes' theorem are appropriate options as they can be formulated in a coherent and consistent manner (*Stuart*, 2010; *Vrugt*, 2016). It means that each time we run the Bayesian method, it should converge to the same distribution. In the Bayesian framework, quantities of interest are modeled as random variables, whose posterior distribution is proportional to the prior distribution times the likelihood. In most situations, closed-form expressions of the posterior distribution are non-existent, so one has to resort to Monte Carlo simulation methods to obtain numerical approximations. Over the past decades, Markov chain Monte Carlo (MCMC) methods have been widely used to assess uncertainties of hydrologic systems conditioned on measurements of observable state variables (*Shi et al.*, 2014; *Smith & Marshall*, 2008; *Vrugt*, 2016; *Zeng*

*et al.*, 2012; *Zeng et al.*, 2018; *Zhang et al.*, 2016; *Zhang et al.*, 2015). However, MCMC has to sufficiently explore the parameter space to obtain reliable estimation results. This usually requires a large number of model evaluations. While it is common that the execution time of a single forward model simulation can be in hours or even longer, e.g., for distributed groundwater models (*Elshall & Tsai*, 2014; *Keating et al.*, 2010). In this situation, the computational cost of MCMC simulation will be prohibitive, especially for high-dimensional problems.

To alleviate the computational cost, one can use a CPU-efficient low-fidelity model (denoted by $f_L(\mathbf{m})$, where $\mathbf{m}$ signifies the model parameters) of the original model (denoted by $f_H(\mathbf{m})$, i.e., the high-fidelity model in this paper) in the MCMC simulation. A low-fidelity model could be a data-driven surrogate based on interpolation or regression, a numerical model that considers fewer processes or has a lower numerical precision (e.g., with a coarser discretization) or is constructed by projecting high-dimensional variables onto their low-dimensional subspace, etc. (*Asher et al.*, 2015; *Mo et al.*, 2017; *Razavi et al.*, 2012; *Smith*, 2014; *Zeng et al.*, 2018; *Zhang et al.*, 2017). Inevitably, using the low fidelity model $f_L(\mathbf{m})$ can introduce some bias if no error model is considered (*Forrester & Keane*, 2009; *Razavi et al.*, 2012). To address this issue, one popular approach is to further evaluate $f_H(\mathbf{m})$ in a two-stage manner, i.e., first explore the parameter space sufficiently using $f_L(\mathbf{m})$ with a low computational cost at stage one, then use $f_H(\mathbf{m})$ to correctly sample the target distribution at stage two (*Efendiev et al.*, 2005; *Laloy et al.*, 2013; *Zeng et al.*, 2012; *Zhang et al.*, 2015). In this way, many unnecessary evaluations of $f_H(\mathbf{m})$ can be avoided. More recently, another approach that adaptively refines a data-driven surrogate over the posterior distribution is proposed, which has shown to be highly efficient (*Gong & Duan*, 2017; *Ju et al.*, 2018; *Zhang et al.*, 2016).

In the above approaches, nevertheless, the correlation relationship between $f_H(\mathbf{m})$ and $f_L(\mathbf{m})$ is not utilized, which leaves some potential untouched. One method that is suitable for fusing $f_H(\mathbf{m})$ and $f_L(\mathbf{m})$ into an integrated system is the multi-fidelity simulation (*Kennedy & O'hagan*, 2000). Given a small number of $f_H(\mathbf{m})$ evaluations and a much larger number of $f_L(\mathbf{m})$ evaluations, the multi-fidelity simulation can take advantage of both the efficiency of $f_L(\mathbf{m})$ and the accuracy of $f_H(\mathbf{m})$. The integrated system can be constructed with many methods, e.g., polynomial chaos expansion (*Narayan et al.*, 2014; *Palar et al.*, 2016; *Zhu et al.*, 2014) and Gaussian process (GP) (*Kennedy & O'hagan*, 2000; *Le Gratiet & Garnier*, 2014; *Parussini et al.*, 2017; *Raissi et al.*, 2017). GP is a generic supervised learning method that uses a (multivariate) Gaussian distribution to predict the quantity of interest based on a set of training data (*Williams & Rasmussen*, 2006). GP has been widely used in hydrologic science to simulate the input-output relationship of the system model (*Asher et al.*, 2015; *Sun et al.*, 2014; *Zhang et al.*, 2016) and the model structural error (*Xu & Valocchi*, 2015; *Xu et al.*, 2017), etc. As GP is very flexible (i.e., we can specify different covariance functions) and can provide the uncertainty estimation (variance) of the system output, we adopt it in this paper to simulate the input-output relationship of $f_H(\mathbf{m})$ with multi-fidelity simulation. In the multi-fidelity GP framework, the correlation between $f_H(\mathbf{m})$ and $f_L(\mathbf{m})$ is rigorously considered with the covariance

(kernel) functions, and the corresponding hyperparameters are estimated with an optimization method conditioned on simulation data of both $f_H(\mathbf{m})$ and $f_L(\mathbf{m})$.

To improve the performance of the multi-fidelity GP system, people would like to acquire the training data adaptively via active learning. One intuitive strategy is to add a new set of training data that have the largest output variance (*Raissi et al.*, 2017). Other strategies can be more sophisticated, e.g., the expected informativeness of candidate points (*Mackay*, 1992), the mutual information criterion (*Krause et al.*, 2008), and the integrated posterior variance (*Gorodetsky & Marzouk*, 2016), just name a few. How to acquire the optimal new training data adaptively for a specific problem is still an open problem. In many cases, the posterior occupies a very small proportion of the prior region. If the multi-fidelity system is built over the whole prior distribution, its accuracy cannot be guaranteed if only a limited number of $f_H(\mathbf{m})$ evaluations are affordable. Considering that the posterior distribution is the region of interest, it is wise to distribute most of the computational budget (i.e., $f_H(\mathbf{m})$ evaluations) therein. Based on this idea, we propose an adaptive multi-fidelity MCMC (AMF-MCMC) algorithm for efficient inverse modeling of hydrologic systems in this paper. Here, we iteratively run MCMC with a GP system constructed with multi-fidelity simulation to search the posterior region and accordingly add new $f_H(\mathbf{m})$ and $f_L(\mathbf{m})$ data to refine the GP system locally. Gradually, the GP system will be accurate enough in the posterior region. Finally, we can obtain an accurate estimate of the posterior distribution. In the AMF-MCMC algorithm, most of the $f_H(\mathbf{m})$ simulations are run near or within the posterior region, thus very little computational budget is wasted.

In recent years, multi-fidelity simulation methods have become increasingly popular in forward uncertainty quantification problems in hydrologic science (*Linde et al.*, 2017; *Lu et al.*, 2016; *Moslehi et al.*, 2015), while the application of these methods in Bayesian inference is very limited. Here we adopt the method originally developed by *Kennedy & O'hagan* (2000) to build the multi-fidelity GP system. The novelty of this paper is that we propose an efficient framework to build the multi-fidelity GP system adaptively in the posterior distribution, not the traditional way in the prior distribution. The method proposed in this work can accurately sample the posterior distribution, but only requires a small number of $f_H(\mathbf{m})$ evaluations. To our best knowledge, this framework (i.e., refine the multi-fidelity surrogate adaptively over the posterior distribution) is rather new and it has value in both theoretical and practical aspects.

The remainder of this paper is organized as follows. In Section 2, we provide a detailed formulation of the AMF-MCMC algorithm. Then, its performance is demonstrated by three numerical case studies in Section 3. Finally, some conclusions and discussions are provided in Section 4.

## 2. Methods

### 2.1. Bayesian inference with DREAM$_{(ZS)}$

For simplicity, here we represent the observation process of the hydrologic system of concern in the following compact form:

$$\widetilde{\mathbf{d}} = f(\mathbf{m}) + \boldsymbol{\varepsilon}, \qquad (1)$$

where $\widetilde{\mathbf{d}}$ is a $N_d$-vector for the measurements, $f(\cdot)$ is the system model, $\mathbf{m}$ is a

$N_\text{m}$-vector for the unknown model parameters, and $\boldsymbol{\varepsilon}$ is a $N_\text{d}$-vector for the measurement errors. Before obtaining the measurements $\widetilde{\mathbf{d}}$, our knowledge about the unknown model parameters $\mathbf{m}$ is represented by the prior distribution, $p(\mathbf{m})$. When $\widetilde{\mathbf{d}}$ is available, we can update our knowledge about $\mathbf{m}$ with the posterior distribution, $p(\mathbf{m}|\widetilde{\mathbf{d}})$, according to Bayes' theorem:

$$p(\mathbf{m}|\widetilde{\mathbf{d}}) = \frac{p(\mathbf{m})p(\widetilde{\mathbf{d}}|\mathbf{m})}{p(\widetilde{\mathbf{d}})}, \qquad (2)$$

where $L(\mathbf{m}|\widetilde{\mathbf{d}}) \equiv p(\widetilde{\mathbf{d}}|\mathbf{m})$ is the likelihood function that quantifies the mismatch between the model outputs $f(\mathbf{m})$ and the measurements $\widetilde{\mathbf{d}}$, $p(\widetilde{\mathbf{d}}) = \int p(\widetilde{\mathbf{d}}|\mathbf{m})p(\mathbf{m})d\mathbf{m}$ is the evidence. When the measurement errors are assumed to fit multivariate Gaussian distribution, the likelihood $L(\mathbf{m}|\widetilde{\mathbf{d}})$ can be expressed as:

$$L(\mathbf{m}|\widetilde{\mathbf{d}}) = \frac{1}{(2\pi)^{N_\text{d}/2}|\boldsymbol{\Sigma}|^{1/2}} \exp\left\{-\frac{1}{2}[\widetilde{\mathbf{d}} - f(\mathbf{m})]^\text{T}\boldsymbol{\Sigma}^{-1}[\widetilde{\mathbf{d}} - f(\mathbf{m})]\right\}, \qquad (3)$$

where $\boldsymbol{\Sigma}$ is the covariance of the measurement errors, $|\cdot|$ signifies the determinant operator.

In most situations, analytical forms of $p(\mathbf{m}|\widetilde{\mathbf{d}})$ do not exist. Here we resort to an efficient MCMC algorithm, i.e., DREAM$_\text{(ZS)}$ (*Laloy & Vrugt*, 2012; *Laloy et al.*, 2013; *Vrugt*, 2016), to explore the parameter space and estimate $p(\mathbf{m}|\widetilde{\mathbf{d}})$ numerically. MCMC works by constructing a Markov chain that gradually converges to the posterior distribution. To better explore the parameter space, $N_\text{c}$ parallel chains are generated and evolved simultaneously by DREAM$_\text{(ZS)}$. At iteration $t$, based on an archive of thinned chain history, $\mathbf{Z}$, the parallel direction jump and the snooker jump are used to update the previous state in the $i$th chain (i.e., $\mathbf{m}_{t-1}^i$) to obtain the proposal state, $\mathbf{m}_\text{p}^i$.

Then, we will compare the acceptance rate, $p_{\text{acc}} = \min[1, p(\mathbf{m}_p^i|\tilde{\mathbf{d}})/p(\mathbf{m}_{t-1}^i|\tilde{\mathbf{d}})]$, with a random sample $u$ drawn from the uniform distribution, $\mathcal{U}(0,1)$. If $p_{\text{acc}} > u$, we will accept $\mathbf{m}_p^i$ and let $\mathbf{m}_t^i = \mathbf{m}_p^i$; Otherwise, we will reject $\mathbf{m}_p^i$ and let $\mathbf{m}_t^i = \mathbf{m}_{t-1}^i$. Every $T_{\text{thin}}$ iterations, we will append the current $N_c$ states in the Markov chains to the archive $\mathbf{Z}$. When the Markov chains converge to the stationary regime, we can view states in the chains as random samples drawn from the posterior distribution. For more details about DREAM$_{(ZS)}$ and related algorithms, one can refer to (*Vrugt*, 2016).

**2.2. Multi-fidelity simulation with Gaussian process**

To obtain reliable estimations, MCMC needs to sufficiently explore the parameter space, which generally requires a large number of model evaluations. When $f(\mathbf{m})$ is CPU-intensive (then $f(\mathbf{m})$ is called the high-fidelity model and represented by $f_H(\mathbf{m})$ thereafter), the computational cost of MCMC simulation will be prohibitive. In this situation, a CPU-efficient low-fidelity model $f_L(\mathbf{m})$ is usually adopted. To balance accuracy and efficiency, it is desirable to construct an integrated system for the MCMC simulation through fusing the information provided by a small number of $f_H(\mathbf{m})$ evaluations and a much larger number of $f_L(\mathbf{m})$ evaluations. This can be realized with the auto-regressive model (*Kennedy & O'hagan*, 2000):

$$u_H(\mathbf{m}) = \rho u_L(\mathbf{m}) + \delta(\mathbf{m}), \tag{4}$$

where $u_L(\mathbf{m}) \sim \mathcal{GP}(0, k_1(\mathbf{m}, \mathbf{m}'; \boldsymbol{\phi}_1))$ and $\delta(\mathbf{m}) \sim \mathcal{GP}(0, k_2(\mathbf{m}, \mathbf{m}'; \boldsymbol{\phi}_2))$ are two independent GPs; $k_i(\mathbf{m}, \mathbf{m}'; \boldsymbol{\phi}_i)$ are the covariance functions with hyperparameters $\boldsymbol{\phi}_i$, for $i = 1,2$; $\rho$ is the cross-correlation coefficient; $\mathbf{m}$ and $\mathbf{m}'$ are two arbitrary

samples in the parameter space. Then $u_H(\mathbf{m})$ has the following form (*Raissi et al.*, 2017):

$$u_H(\mathbf{m}) \sim \mathcal{GP}(0, k(\mathbf{m}, \mathbf{m}'; \boldsymbol{\phi})), \tag{5}$$

where $k(\mathbf{m}, \mathbf{m}'; \boldsymbol{\phi}) = \rho^2 k_1(\mathbf{m}, \mathbf{m}'; \boldsymbol{\phi}_1) + k_2(\mathbf{m}, \mathbf{m}'; \boldsymbol{\phi}_2)$, and $\boldsymbol{\phi} = [\boldsymbol{\phi}_1, \boldsymbol{\phi}_2, \rho]$. Consequently, we have:

$$\begin{bmatrix} u_L(\mathbf{m}) \\ u_H(\mathbf{m}) \end{bmatrix} \sim \mathcal{GP}\left(0, \begin{bmatrix} k_{LL} & k_{LH} \\ k_{HL} & k_{HH} \end{bmatrix}\right), \tag{6}$$

where $k_{LL} = k_1(\mathbf{m}, \mathbf{m}'; \boldsymbol{\phi}_1)$, $k_{LH} = k_{HL}^T = \rho k_1(\mathbf{m}, \mathbf{m}'; \boldsymbol{\phi}_1)$, and $k_{HH} = k(\mathbf{m}, \mathbf{m}'; \boldsymbol{\phi})$. In this paper, we adopt the commonly used squared exponential covariance function (*Williams & Rasmussen*, 2006) for both $k_1$ and $k_2$:

$$k_i(\mathbf{m}, \mathbf{m}'; \boldsymbol{\phi}_i) = \sigma_i^2 \exp\left[-\frac{1}{2} \sum_{n=1}^{N_m} \frac{(\mathbf{m}_n - \mathbf{m}'_n)^2}{l_{n,i}^2}\right], \tag{7}$$

where $\boldsymbol{\phi}_i = [\sigma_i^2, (l_{n,i}^2)_{n=1}^{N_m}]$ are the hyperparameters of the covariance functions $k_i$, for $i = 1,2$.

When we have $N_L$ sets of training data obtained by simulating $f_L(\mathbf{m})$ at $N_L$ parameter samples and $N_H$ sets of training data obtained by simulating $f_H(\mathbf{m})$ at $N_H$ parameter samples (here $N_L$ and $N_H$ are two positive integers), we can estimate the hyperparameters $\boldsymbol{\phi}$ by minimizing the negative log marginal likelihood:

$$NL = -\log p(\mathbf{D}|\mathbf{M}, \boldsymbol{\phi}) = \frac{1}{2} \mathbf{D}^T \mathbf{K}^{-1} \mathbf{D} + \frac{1}{2} \log|\mathbf{K}| + \frac{N_L + N_H}{2} \log(2\pi), \tag{8}$$

where

$$\mathbf{K} = \begin{bmatrix} k_{LL}(\mathbf{M}_L, \mathbf{M}_L) + \sigma_L^2 \mathbf{I}_{N_L} & k_{LH}(\mathbf{M}_L, \mathbf{M}_H) \\ k_{HL}(\mathbf{M}_H, \mathbf{M}_L) & k_{HH}(\mathbf{M}_H, \mathbf{M}_H) + \sigma_H^2 \mathbf{I}_{N_H} \end{bmatrix}; \tag{9}$$

$\mathbf{M}_L = [\mathbf{m}_1, \dots, \mathbf{m}_{N_L}]$ are $N_L$ parameter samples for $f_L(\mathbf{m})$, $\mathbf{D}_L =$

$[f_L(\mathbf{m}_1), \ldots, f_L(\mathbf{m}_{N_L})]$ are the corresponding low-fidelity model outputs; $\mathbf{M}_H = [\mathbf{m}_1, \ldots, \mathbf{m}_{N_H}]$ are $N_H$ parameter samples for $f_H(\mathbf{m})$, $\mathbf{D}_H = [f_H(\mathbf{m}_1), \ldots, f_H(\mathbf{m}_{N_H})]$ are the corresponding high-fidelity model outputs; $\mathbf{M} = [\mathbf{M}_L\ \mathbf{M}_H]$, $\mathbf{D} = [\mathbf{D}_L\ \mathbf{D}_H]$; $\mathbf{I}_{N_L}$ and $\mathbf{I}_{N_H}$ are identity matrices of size $N_L$ and $N_H$, respectively. Here the hyperparameters to be estimated are $[\boldsymbol{\phi}, \sigma_L, \sigma_H]$. The above equation (8) calculates the goodness-of-fit given the training data and a set of hyperparameters, and it can be easily derived from the (logarithmic) probability density function of multivariate Gaussian distribution (it is noted here that the most general zero-mean function of GP is used) (*Williams & Rasmussen*, 2006).

Minimizing equation (8) defines a non-convex optimization problem. In practice, we use the trust-region algorithm, which is based on the interior-reflective Newton method described in (*Coleman & Li*, 1996), to find the minimum solution. Although it does not guarantee convergence to a global optimum, it is usually sufficient for obtaining a good solution. The most computationally intensive part of solving the optimization problem is associated with inverting the covariance matrix $\mathbf{K}$, which scales cubically with the number of training data. This is a well-known limitation of GP, but it has been effectively addressed with techniques like sparse GP (*Quinonero-Candela & Rasmussen*, 2005; *Snelson & Ghahramani*, 2006). As will be demonstrated in Section 2.3, the method proposed in this paper does not require a large number of training data. Moreover, the process of training GPs can be accelerated by adopting parallel computation. Thus, the computational cost of training GPs will not be a big problem.

After training the GP system conditioned on the multi-fidelity data $\mathbf{D}$ (i.e., after obtaining the optimal hyperparameters by minimizing the objective function defined in equation (8)), we can obtain a new GP system represented by $\tilde{g}(\mathbf{m})$, which can be used to predict the model output given an arbitrary parameter set $\mathbf{m}^*$:

$$\tilde{g}(\mathbf{m}^*) = u_\mathrm{H}(\mathbf{m}^*)|\mathbf{D} \sim \mathcal{N}(\mu_\mathrm{gp}, \sigma_\mathrm{gp}^2), \qquad (10)$$

where $\mu_\mathrm{gp} = \mathbf{a}\mathbf{K}^{-1}\mathbf{D}$ is the mean estimate, $\sigma_\mathrm{gp}^2 = k_\mathrm{HH}(\mathbf{m}^*, \mathbf{m}^*) - \mathbf{a}\mathbf{K}^{-1}\mathbf{a}^\mathrm{T}$ is the estimation variance, and $\mathbf{a} = [k_\mathrm{HL}(\mathbf{m}^*, \mathbf{M}_\mathrm{L})\ k_\mathrm{HH}(\mathbf{m}^*, \mathbf{M}_\mathrm{H})]$. Then the multi-fidelity system $\tilde{g}(\mathbf{m})$ can be used in MCMC simulations to gain efficiency.

[Figure 1]

As the multi-fidelity GP (MF-GP) can fuse the efficiency of $f_\mathrm{L}(\mathbf{m})$ and the accuracy of $f_\mathrm{H}(\mathbf{m})$, it is advantageous over the GP system constructed only on data from $f_\mathrm{H}(\mathbf{m})$ (here we call it the single-fidelity GP, i.e., SF-GP). As shown in Figure 1, in a one-dimensional problem, with three training data from the high-fidelity model and twenty training data from the low-fidelity model, we can build an MF-GP system that is more accurate than the SF-GP system based on four training data from the high-fidelity model, especially in the area that is far away from the high-fidelity training data. Here the high-fidelity model is $f_\mathrm{H}(m) = \sin(m)$, and the low-fidelity model is obtained by adding an error term, $-0.1m - 0.1$, to the high-fidelity model, i.e., $f_\mathrm{L}(m) = \sin(m) - 0.1m - 0.1$. In this simple case, both the high- and low-fidelity models are very quick to simulate, and this case is only used to demonstrate the performance of multi-fidelity GP.

In the above approach, we only utilize data from two levels of model fidelity. If $s$-levels of data $\{\mathbf{D}_t(\mathbf{M})\}_{t=1}^{s}$ sorted by increasing fidelity are available, we can readily extend the auto-regressive scheme:

$$u_t(\mathbf{m}) = \rho_{t-1} u_{t-1}(\mathbf{m}) + \delta_t(\mathbf{m}), t = 2, \ldots, s, \tag{11}$$

where $\delta_t(\mathbf{m})$ is a Gaussian process independent of $\{u_{t-1}(\mathbf{m}), \ldots, u_1(\mathbf{m})\}$, $\rho_{t-1}$ is the cross-correlation coefficient. Moreover, we can consider complex, nonlinear relationship between $u_{t-1}(\mathbf{m})$ and $u_t(\mathbf{m})$ (*Perdikaris et al.*, 2017):

$$u_t(\mathbf{m}) = q_{t-1}(u_{t-1}(\mathbf{m})) + \delta_t(\mathbf{m}), t = 2, \ldots, s, \tag{12}$$

where $q_{t-1}(\cdot)$ is a nonlinear mapping.

It is noted here that the autoregressive model considers the relationship between the high- and low-fidelity models, not the input-output relationship of the hydrologic system. In practice, the low-fidelity model should resemble and capture the right trend of the high-fidelity model, otherwise the low-fidelity model would be useless. In that case, using a linear autoregressive model is still applicable even in nonlinear problems. For the sake of clarity, the nonlinear mapping described above is not adopted in the present paper. For more details about GP and its construction with multi-fidelity data, one can refer to (*Kennedy & O'hagan*, 2000; *Parussini et al.*, 2017; *Raissi et al.*, 2017; *Williams & Rasmussen*, 2006).

### 2.3. The adaptive multi-fidelity MCMC algorithm

Generally, we can only afford a limited number of $f_H(\mathbf{m})$ evaluations. If the multi-fidelity system $\tilde{g}(\mathbf{m})$ is constructed over the whole prior distribution, its accuracy cannot be guaranteed. In MCMC simulations, our concern is the posterior distribution. Thus, it is crucial that $\tilde{g}(\mathbf{m})$ is accurate enough therein, but there is no need to ensure its accuracy elsewhere. Below we propose an adaptive multi-fidelity MCMC (AMF-MCMC) algorithm which adaptively refines $\tilde{g}(\mathbf{m})$ over the posterior distribution and finally obtains an accurate estimate of $p(\mathbf{m}|\tilde{\mathbf{d}})$.

The AMF-MCMC algorithm first builds an initial multi-fidelity system $\tilde{g}_0(\mathbf{m})$ conditioned on $N_\text{L}$ evaluations of $f_\text{L}(\mathbf{m})$ and $N_\text{H}$ evaluations of $f_\text{H}(\mathbf{m})$, where $N_\text{H}$ is usually a small integer and $N_\text{L}$ is much larger than $N_\text{H}$. With $\tilde{g}_0(\mathbf{m})$, we can run DREAM$_\text{(ZS)}$ to sufficiently explore the parameter space, which can be done very quickly. Here the variance of the GP system ($\sigma_\text{gp}^2$) is considered in the MCMC simulation by augmenting it with the variance of the measurement error ($\sigma_\text{meas}^2$), i.e., $\sigma_\text{total}^2 = \sigma_\text{gp}^2 + \sigma_\text{meas}^2$, and using $\sigma_\text{total}^2$ in the likelihood function, $L(\mathbf{m}|\tilde{\mathbf{d}})$. Then we can draw two random samples, $\mathbf{m}_\text{H}^\text{p}$ and $\mathbf{m}_\text{L}^\text{p}$, from the approximated posterior, $\tilde{p}_0(\mathbf{m}|\tilde{\mathbf{d}})$, which are expected to be much closer to the posterior region than the prior samples. By utilizing the new data $f_\text{H}(\mathbf{m}_\text{H}^\text{p})$ and $f_\text{L}(\mathbf{m}_\text{L}^\text{p})$, the multi-fidelity GP system can be refined locally. This process will be further iterated $I_\text{max} - 1$ times. Since it is very cheap to evaluate $f_\text{L}(\mathbf{m})$, we can also collect more samples from the approximated posterior for $f_\text{L}(\mathbf{m})$ simulations in GP refinement at each iteration. However, an increased number of training data will inevitably increase the computational cost in constructing the multi-fidelity GP system. To address this issue, one practical strategy is to remove some "bad" data according to the closeness of the training data to the measurements (e.g., the likelihood function defined in equation (3)).

In the AMF-MCMC algorithm, we gradually improve the accuracy of the multi-fidelity GP system by iteratively adding new training data from $f_H(\mathbf{m})$ and $f_L(\mathbf{m})$ simulations in the posterior distribution. At each iteration, the new training data are obtained by running MCMC based on the previous GP system to sufficiently explore the parameter space. As the GP error is considered, the MCMC simulation results will not be overconfident and biased (although the posterior will be wider). Based on both the new and old training data, we can obtain an updated GP system that is slightly more accurate in the posterior distribution. Finally, we can obtain a GP system that is locally accurate (i.e., with a negligible $\sigma_{gp}^2$) in the posterior region and an accurate estimate of $p(\mathbf{m}|\widetilde{\mathbf{d}})$.

After each DREAM$_{(ZS)}$ simulation, the thinned chain history $\mathbf{Z}$ will be saved and used in the next DREAM$_{(ZS)}$ simulation. Actually, this treatment tailors the $I_{\max}+1$ $\tilde{g}(\mathbf{m})$-based MCMC simulations into an integrated one, which is beneficial to better explore the parameter space. The complete scheme of the AMF-MCMC algorithm is given in Algorithm 1.

---

**Algorithm 1** The adaptive multi-fidelity MCMC algorithm.

1. Draw $N_L$ random samples from the prior distribution, $\mathbf{M}_L = [\mathbf{m}_1, \dots, \mathbf{m}_{N_L}]$, calculate $\mathbf{D}_L = [f_L(\mathbf{m}_1), \dots, f_L(\mathbf{m}_{N_L})]$.
2. Draw $N_H$ random samples from the prior distribution, $\mathbf{M}_H = [\mathbf{m}_1, \dots, \mathbf{m}_{N_H}]$, calculate $\mathbf{D}_H = [f_H(\mathbf{m}_1), \dots, f_H(\mathbf{m}_{N_H})]$. Here $N_H \ll N_L$.
3. Build the initial GP system $\tilde{g}_0(\mathbf{m})$ with multi-fidelity simulation conditioned on $[\mathbf{M}_L \ \mathbf{M}_H]$ and $[\mathbf{D}_L \ \mathbf{D}_H]$.
4. Run MCMC with $\tilde{g}_0(\mathbf{m})$, obtain $\tilde{p}_0(\mathbf{m}|\widetilde{\mathbf{d}})$.
5. **for** $i = 1, \dots, I_{\max}$ **do**
    Draw two random samples, $\mathbf{m}_H^p$ and $\mathbf{m}_L^p$, from $\tilde{p}_{i-1}(\mathbf{m}|\widetilde{\mathbf{d}})$, let $\mathbf{M}_H = [\mathbf{M}_H \ \mathbf{m}_H^p]$, $\mathbf{D}_H = [\mathbf{D}_H \ f_H(\mathbf{m}_H^p)]$, $\mathbf{M}_L = [\mathbf{M}_L \ \mathbf{m}_L^p]$, and $\mathbf{D}_L = [\mathbf{D}_L \ f_L(\mathbf{m}_L^p)]$.
    Build the GP system $\tilde{g}_i(\mathbf{m})$ conditioned on $[\mathbf{M}_L \ \mathbf{M}_H]$ and $[\mathbf{D}_L \ \mathbf{D}_H]$.
    Run MCMC with $\tilde{g}_i(\mathbf{m})$ and previous MCMC simulation results, obtain $\tilde{p}_i(\mathbf{m}|\widetilde{\mathbf{d}})$.

    **end for**

6. The posterior is approximated with $\tilde{p}_{I_{\max}}(\mathbf{m}|\widetilde{\mathbf{d}})$.

## 3. Illustrative examples

### 3.1. Example 1: Estimation of soil hydraulic and thermal parameters

In example 1, we first test the performance of the AMF-MCMC algorithm in estimating soil hydraulic and thermal parameters in a single ring infiltration experiment (*Nakhaei & Šimůnek*, 2014). Here the processes of water flow and heat transport are considered. As shown in Figure 2, the flow domain is $100\text{cm} \times 200\text{cm}$. The initial conditions for water content and temperature in the domain are $0.100\text{cm}^3\,\text{cm}^{-3}$ and $17.5\,°C$, respectively. The domain has three types of boundary conditions, i.e., impervious condition at the two lateral boundaries and part of the upper boundary (represented by the green lines in Figure 2), free drainage condition at the lower boundary (represented by the blue line in Figure 2) and constant temperature (61.0°C) and water content ($0.430\,\text{cm}^3\,\text{cm}^{-3}$) conditions at part of the upper boundary (represented by the red line in Figure 2).

[Figure 2]

With the initial and boundary conditions prescribed above, we can simulate unsaturated water flow in the domain by numerically solving the Richards' equation with HYDRUS-2D (*Šimůnek et al.*, 2008):

$$\frac{\partial \theta}{\partial t} = \frac{\partial}{\partial x}\left[K(h)\frac{\partial h}{\partial x}\right] + \frac{\partial}{\partial z}\left[K(h)\frac{\partial h}{\partial z} + K(h)\right], \tag{13}$$

where $\theta[L^3L^{-3}]$ is volumetric water content of the soil; $t[T]$ is time; $x[L]$ and $z[L]$ are distances along the horizontal and vertical directions; $h[L]$ is pressure head; $K(h)[LT^{-1}]$ is hydraulic conductivity, which is a function of $h$ (*Mualem*, 1976; *Van*

*Genuchten*, 1980):

$$K(h) = K_s S_e^l \left[1 - \left(1 - S_e^{1/m}\right)^m\right]^2, \quad (14)$$

where $K_s [LT^{-1}]$ is saturated hydraulic conductivity, $S_e[-]$ is effective saturation:

$$S_e = \frac{\theta - \theta_r}{\theta_s - \theta_r} = \begin{cases} \frac{1}{(1 + |\alpha h|^n)^m} & h < 0 \\ 1 & h \geq 0 \end{cases}, \quad (15)$$

where $\theta_r$ and $\theta_s$ are residual and saturated water content $[L^3 L^{-3}]$; $l[-]$ is a pore-connectivity parameter; $\alpha[L^{-1}]$, $n[-]$ and $m = (1 - 1/n)[-]$ are empirical shape parameters, respectively.

Based on the simulation results of unsaturated water flow, we can further simulate heat transport by numerically solving the following governing equation (*Sophocleous*, 1979) with HYDRUS-2D:

$$C(\theta) \frac{\partial T}{\partial t} = \frac{\partial}{\partial z}\left[\lambda_{xz}(\theta) \frac{\partial T}{\partial x}\right] - C_w q_z \frac{\partial T}{\partial z}, \quad (16)$$

where $C(\theta) = C_n \theta_n + C_o \theta_o + C_w \theta$ is volumetric heat capacity of soil $[ML^{-1}T^{-2}K^{-1}]$, $C_n$, $C_o$ and $C_w$ are volumetric heat capacities of solid phase, organic phase and liquid phase $[ML^{-1}T^{-2}K^{-1}]$, $\theta_n$ and $\theta_o$ are fraction of solid phase and organic phase $[L^3 L^{-3}]$, respectively; $T[K]$ is temperature; $\lambda_{xz}(\theta)[MLT^{-3}K^{-1}]$ is apparent thermal conductivity:

$$\lambda_{xz}(\theta) = \lambda_x C_w |q| \delta_{xz} + (\lambda_z - \lambda_x) C_w \frac{q_x q_z}{|q|} + \lambda_0(\theta) \delta_{xz}, \quad (17)$$

where $\lambda_x$ and $\lambda_z$ are components of thermal dispersivity in the horizontal and vertical directions $[L]$; $q[LT^{-1}]$ is fluid flux density with absolute value $|q|$, and components in the horizontal and vertical directions, $q_x$ and $q_z$, respectively; $\delta_{xz}$ is Kronecker delta function; $\lambda_0(\theta) = b_1 + b_2 \theta + b_3 \theta^{0.5}$ is thermal conductivity of soil in the absence of flow, where $b_1$, $b_2$ and $b_3$ are empirical parameters, $[\text{MLT}^{-3}\text{K}^{-1}]$. In this example, the total simulation time is 10 hours. For more details about the model descriptions, one can refer to (*Nakhaei & Šimůnek*, 2014).

[Table 1]

Here the unknown model parameters are $\alpha$, $n$, $K_s$, $b_1$, $b_2$ and $b_3$, which are assumed to be homogeneous in the domain and fit multivariate uniform prior distribution (Table 1). While other parameters are assumed to be known, as listed in Table 2. To infer the six unknown model parameters, we collect pressure head and temperature measurements at $t = 1, 2, \ldots 10$ hour at three locations denoted by the magenta circles in Figure 2. The measurements are generated from one set of true model parameters (Table 1) evaluated with $f_\text{H}(\mathbf{m})$ and perturbed with additive measurement errors, $\boldsymbol{\varepsilon} \sim \mathcal{N}(\mathbf{0}, \boldsymbol{\sigma}^2)$. Here the standard deviation of the measurement errors for pressure head is $\sigma_h = 1\text{cm}$, the standard deviation of the measurement errors for temperature is $\sigma_T = 0.5°\text{C}$. Then we use the Gaussian likelihood function defined in equation (3) to evaluate the goodness-of-fit between the model simulation and measurement data.

[Table 2]

[Figure 3]

In this example, the high- and low-fidelity models are constructed with different levels of discretization. In $f_\text{H}(\mathbf{m})$, the flow domain is evenly discretized into $41 \times 41$ grids. While in $f_\text{L}(\mathbf{m})$, there are $21 \times 21$ grids. In a Dell Precision T7610 workstation (Intel Xeon Processor E5-2680 v2 @ 2.80GHz; 512GB RAM; Windows 10 Pro., 64

Bit), the average time of evaluating $f_H(\mathbf{m})$ once is about 54 seconds, while the average time of evaluating $f_L(\mathbf{m})$ once is about 5.6 seconds. We admit here that, the high-fidelity model used in this test case is not actually very time-consuming, so that we can implement the high-fidelity MCMC simulation to obtain the reference results to verify the performance of the AMF-MCMC algorithm. As shown in Figure 3, although $f_L(\mathbf{m})$ can capture the main trend of $f_H(\mathbf{m})$, systematic errors do exist.

[Figure 4]

We first run MCMC with $f_H(\mathbf{m})$ and $f_L(\mathbf{m})$ respectively to approximate the posterior. Here the two MCMC simulations have 4 parallel chains with 4,000 iterations, which means 16,000 model evaluations for each simulation. Using the $\hat{R}$-statistic proposed in (*Brooks & Gelman*, 1998; *Gelman & Rubin*, 1992), we can monitor the convergence of the $f_H(\mathbf{m})$- and $f_L(\mathbf{m})$-MCMC simulations. In each panel of Figure 4, traces of the $\hat{R}$-statistic of different parameters are coded with different colors. The red dashed line in both panels demarcates the threshold of 1.2 below which the chains are assumed to have converged to a stationary distribution. As shown in Figure 5, the $f_H(\mathbf{m})$-MCMC simulation can identify the unknown model parameters properly, while the estimation results of the $f_L(\mathbf{m})$-MCMC simulation are significantly biased, especially for the first three parameters. Thus, although using a low-fidelity model in the MCMC simulation can gain computational efficiency, the estimation accuracy cannot be guaranteed.

[Figure 5]

We then run the AMF-MCMC algorithm with the same set of measurements to approximate the posterior. Here the initial number of $f_L(\mathbf{m})$ evaluations is $N_L = 200$, and the initial number of $f_H(\mathbf{m})$ evaluations is $N_H = 30$. These initial samples are randomly drawn from the prior distribution. Conditioned on these multi-fidelity data, we build the initial GP system, $\tilde{g}_0(\mathbf{m})$, based on which the MCMC simulation (4 parallel chains with 4,000 iterations) can be implemented very quickly. From the approximated posterior, $\tilde{p}_0(\mathbf{m}|\tilde{\mathbf{d}})$, we can draw two random parameter samples, $\mathbf{m}_H^p$ and $\mathbf{m}_L^p$, which are expected to be much closer to the posterior region than the prior samples. Then we can add $f_H(\mathbf{m}_H^p)$ and $f_L(\mathbf{m}_L^p)$ to the existing multi-fidelity training data to refine the GP system locally. The coupled process of GP-based MCMC simulation and GP system refinement is further repeated 69 times (i.e., $I_{\max} = 70$). In Figure 6, we plot the evolution of the variance of the multi-fidelity GP output (i.e., $\hat{\sigma}_{gp}(h)$ and $\hat{\sigma}_{gp}(T)$) and the RMSE between the multi-fidelity GP output and the high-fidelity output (i.e., RMSE($h$) and RMSE($T$)) averaged over 400 posterior samples. It is found that for both hydraulic head and temperature outputs, the accuracy of the multi-fidelity GP system will improve with the iteration. Finally, the multi-fidelity GP system will be rather accurate in the posterior region.

[Figure 6]

As shown in Figure 7, the successively added new parameter samples (blue dots) for $f_H(\mathbf{m})$ evaluations gradually approach to the true values (black crosses), which are the basis of a locally accurate GP system. Based on this GP system, we can obtain a rather accurate approximation of $p(\mathbf{m}|\tilde{\mathbf{d}})$. As shown in Figure 8, the marginal posterior probability density functions (PPDFs) obtained by the $f_H(\mathbf{m})$-MCMC

algorithm and the AMF-MCMC algorithm are almost identical, which confirms the accuracy of the AMF-MCMC algorithm.

[Figure 7]

[Figure 8]

In the AMF-MCMC algorithm, the total numbers of $f_H(\mathbf{m})$ and $f_L(\mathbf{m})$ evaluations are 100 and 270, respectively. With respect to the number of model evaluations, the AMF-MCMC algorithm is much more efficient than both the $f_H(\mathbf{m})$ and $f_L(\mathbf{m})$-based MCMC simulations. It is noted here that, when building a GP system, the training data of model parameters are a matrix of $N_m \times N_{tr}$, while the training data of model outputs should be a vector of $1 \times N_{tr}$, where $N_{tr}$ is the total number of training data. It means that when the model outputs are $N_d$-dimensional, we have to build $N_d$ GP systems for each of the $N_d$ model outputs separately. When we run the high- or low-fidelity model to acquire the training data, we can obtain the $N_d$ model outputs at the same time. So the total number of function evaluations that are used to build GPs for all observation data points is still $N_{tr}$. Nevertheless, the time needed by the multi-fidelity GP system constructions and GP-based MCMC simulations should not be neglected, especially when the number of model outputs $N_d$ is large. To improve the efficiency of the AMF-MCMC algorithm, we can build the $N_d$ GPs in parallel. In our simulations, there are 20 cores available, which can be utilized to greatly accelerate the simulation of the AMF-MCMC algorithm.

[Figure 9]

In the above simulation, we only add one set of $f_H(\mathbf{m})$ data and one set of $f_L(\mathbf{m})$

data at each iteration. As stated in Section 2.3, we can also add more than one set of new data once. In Figure 9, we compare the finally obtained variance of the multi-fidelity GP output (i.e., $\hat{\sigma}_{\text{gp}}(h)$ and $\hat{\sigma}_{\text{gp}}(T)$) and RMSE between the multi-fidelity GP output and the high-fidelity output (i.e., $\text{RMSE}(h)$ and $\text{RMSE}(T)$) when adding different number of new training data once. Here the iteration numbers are set as 70, 35, 14, 7, 5 and 4 for adding 1, 2, 5, 10, 15 and 20 sets of $f_{\text{H}}(\mathbf{m})$ and $f_{\text{L}}(\mathbf{m})$ training data at each iteration, respectively, to make sure that the total numbers of added data sets are roughly equal. It is clear that in this case, adding one set of $f_{\text{H}}(\mathbf{m})$ data and one set of $f_{\text{L}}(\mathbf{m})$ data once can bring about a more accurate multi-fidelity system in the posterior region.

[Figure 10]

In our previous paper (*Zhang et al.*, 2016), we proposed an efficient method to estimate hydrologic model parameters by combing MCMC simulations with a GP surrogate adaptively refined over the posterior distribution. Here we call that method the adaptive Gaussian process MCMC (AGP-MCMC) algorithm. It is obvious that the AMF-MCMC algorithm proposed in the present paper is an extension of the AGP-MCMC algorithm through introducing the adaptive multi-fidelity simulation. As the AMF-MCMC algorithm takes advantage of the correlation relationship between $f_{\text{H}}(\mathbf{m})$ and $f_{\text{L}}(\mathbf{m})$, it generally requires less $f_{\text{H}}(\mathbf{m})$ evaluations than the AGP-MCMC algorithm. In other words, with similar computational cost, the multi-fidelity GP system helps to better explore the parameter space compared to the single-fidelity GP system. To demonstrate this point, we further compare the performances of the two algorithms in this example when limited numbers of $f_H(\mathbf{m})$ and $f_L(\mathbf{m})$ evaluations are affordable. In the AMF-MCMC algorithm, the initial number of $f_L(\mathbf{m})$ evaluations is $N_L = 60$, the initial number of $f_H(\mathbf{m})$ evaluations is $N_H = 10$, and an extra number of $I_{max} = 40$ $f_H(\mathbf{m})$ and $f_L(\mathbf{m})$ evaluations are called to adaptively refine the multi-fidelity GP system in the posterior region. In the AGP-MCMC algorithm, we first evaluate $N_H = 20$ prior samples with $f_H(\mathbf{m})$. Then we further evaluate $I_{max} = 40$ more parameter samples with $f_H(\mathbf{m})$ that are gradually approaching to the posterior region. Since the simulation time of 100 $f_L(\mathbf{m})$ evaluations is approximately equal to the simulation time of 10 $f_H(\mathbf{m})$ evaluations, the computational costs in model evaluations for AMF-MCMC and AGP-MCMC are roughly equal. As shown in Figure 10 (a-c), the adaptively added parameter samples by the AMF-MCMC algorithm converge faster to the true values. For the AGP-MCMC algorithm, only 60 $f_H(\mathbf{m})$ evaluations in total are still not enough, especially for the parameters $b_1$ and $b_3$.

[Figure 11]

As shown in Figure 11, the finally estimated PPDFs of the unknown model parameters by the AMF-MCMC algorithm (blue dashed curves) are overall much closer to the reference results (red curves) than the AGP-MCMC algorithm (magenta dash-dotted curves). Therefore, it is important to leverage the multi-fidelity system based on a large number of $f_L(\mathbf{m})$ simulations to sufficiently explore the parameter space especially at the early stage of the GP-based MCMC simulations. In practical applications, to guarantee the estimation accuracy of the AMF-MCMC algorithm, more evaluations of $f_H(\mathbf{m})$ are suggested. For more details about the AGP-MCMC algorithm, one can refer to (*Zhang et al.*, 2016) and a more recent work by other researchers (*Gong & Duan*, 2017).

**3.2. Example 2: Contaminant source identification with multimodal posterior**

In example 2, we further test the performance of the AMF-MCMC algorithm in solving an inverse problem with multimodal posterior distribution. Here we consider the processes of steady-state saturated groundwater flow and contaminant transport. As shown in Figure 12, the flow domain is $20[L] \times 10[L]$. The upper and lower sides are no-flow boundaries, the left ($h = 12[L]$) and right ($h = 11[L]$) sides are constant-head boundaries, respectively. At the initial time, hydraulic heads in the domain are all $11[L]$ except for the left boundary ($12[L]$).

[Figure 12]

Given the above initial and boundary conditions, we can obtain the flow field through numerically solving the following governing equations with MODFLOW (*Harbaugh et al.*, 2000):

$$\frac{\partial}{\partial x_i}\left(K_i \frac{\partial h}{\partial x_i}\right) = 0, \tag{18}$$

and

$$v_i = -\frac{K_i}{\theta}\frac{\partial h}{\partial x_i}, \tag{19}$$

where $x_i[L]$, $K_i[LT^{-1}]$ and $v_i[LT^{-1}]$ signify distance, hydraulic conductivity and pore water velocity along the respective Cartesian coordinate axis, for $i = 1, 2$; $h[L]$ is hydraulic head; $\theta[-]$ is porosity of the aquifer. Here $K = 8[LT^{-1}]$ and $\theta = 0.25$ are known beforehand. In the steady-state flow field, some amount of contaminant is released from a point source located somewhere in the red dashed rectangle depicted in Figure 12. The contaminant source is characterized by five parameters, i.e., location, $(x_s, y_s)$[L], source strength measured by mass loading rate, $S_s$[MT$^{-1}$], start time of contaminant release, $t_{on}$[T], and end time of the release, $t_{off}$[T]. Then we can obtain contaminant concentration $C$[ML$^{-3}$] at different times and locations through numerically solving the following advection-dispersion equation with MT3DMS (*Zheng & Wang*, 1999):

$$\frac{\partial(\theta C)}{\partial t} = \frac{\partial}{\partial x_i}\left(\theta D_{ij}\frac{\partial C}{\partial x_j}\right) - \frac{\partial}{\partial x_i}(\theta v_i C) + q_s C_s, \tag{20}$$

where $t$[T] is time; $q_s$[T$^{-1}$] and $C_s$[ML$^{-3}$] are volumetric flow rate per unit volume of the aquifer and concentration of the source, respectively; $D_{ij}$[L$^2$T$^{-1}$] are hydrodynamic dispersion coefficient tensors:

$$\begin{cases} D_{11} = (\alpha_L v_1^2 + \alpha_T v_2^2)/|v|, \\ D_{22} = (\alpha_L v_2^2 + \alpha_T v_1^2)/|v|, \\ D_{12} = D_{21} = (\alpha_L - \alpha_T)v_1 v_2/|v|, \end{cases} \tag{21}$$

where $\alpha_L = 0.3$[L] is longitudinal dispersivity, $\alpha_T = 0.03$[L] is transverse dispersivity, and $|v| = \sqrt{v_1^2 + v_2^2}$ is magnitude of the velocity.

In this example, the unknown model parameters to be estimated are the five contaminant source parameters. Here we assume that our prior knowledge about them are rather limited and thus represented by a multivariate uniform distribution (Table 3). To infer the five unknown parameters, we collect concentration measurements at $t = 6, 8, 10, 12, 14$[T] at a well denoted by the blue square in Figure 12. Here the measurements are generated from one set of true model parameters (Table 3) evaluated with $f_H(\mathbf{m})$ and perturbed with additive measurement errors that fit $\mathcal{N}(0, 0.01^2)$.

[Table 3]

In this example, the high-fidelity model $f_H(\mathbf{m})$ is the coupled numerical model built with MODFLOW and MT3DMS. While the low-fidelity model $f_L(\mathbf{m})$ is built with a data-driven surrogate, i.e., the adaptive sparse grid interpolation method proposed by *Klimke & Wohlmuth* (2005), based on 41 evaluations of $f_H(\mathbf{m})$. We compare the simulation results between $f_H(\mathbf{m})$ and $f_L(\mathbf{m})$ with 50 random parameter samples drawn from the prior distribution and obtain $R^2 = 0.949$. It should be pointed out here that, other low-fidelity models, e.g., a numerical model with a coarser discretization, can also be utilized. Considering that data-driven surrogates are widely used in hydrologic science (*Asher et al.*, 2015; *Razavi et al.*, 2012), and many times in MCMC simulations (*Elsheikh et al.*, 2014; *Laloy et al.*, 2013; *Zeng et al.*, 2012; *Zeng et al.*, 2016), here we test the applicability of data-driven surrogates in the AMF-MCMC algorithm. Then we run MCMC with $f_H(\mathbf{m})$ and $f_L(\mathbf{m})$ respectively to approximate the posterior (6 parallel chains and 3,000 iterations). As shown in Figure 13, the results obtained by the $f_L(\mathbf{m})$-MCMC simulation are significantly biased due to the approximation errors, especially for $x_s$ and $t_{\text{off}}$.

[Figure 13]

With the same set of measurements, we further run the AMF-MCMC algorithm to approximate the posterior. Here the initial number of $f_L(\mathbf{m})$ evaluations is $N_L = 300$, and the initial number of $f_H(\mathbf{m})$ evaluations is $N_H = 30$. Then we successively add another $I_{\max} = 50$ parameter samples for both $f_H(\mathbf{m})$ and $f_L(\mathbf{m})$ evaluations to refine the GP system over the posterior distribution. Finally, we can obtain a rather accurate approximation of $p(\mathbf{m}|\widetilde{\mathbf{d}})$. As shown in Figure 14, the bivariate scatter plots of the posterior samples obtained by the $f_\mathrm{H}(\mathbf{m})$-MCMC algorithm and the AMF-MCMC algorithm are almost identical, which indicates the accuracy of the AMF-MCMC algorithm. Moreover, the bimodality of $y_\mathrm{s}$ is well identified by both algorithms.

[Figure 14]

Here we also test how the parametric uncertainty affects the credible intervals of the model simulation and prediction for the different MCMC approaches. In this case, the quantity of interest (QoI) is the concentration at a well denoted by the red square in Figure 12 at a future time $t = 16[\mathrm{T}]$. To quantify the predictive uncertainty of the QoI, we further run $f_\mathrm{H}(\mathbf{m})$ to the future time given posterior samples obtained by the $f_\mathrm{H}(\mathbf{m})$-MCMC simulation, the $f_\mathrm{L}(\mathbf{m})$-MCMC simulation and the AMF-MCMC simulation, respectively. As shown in Figure 15, the AMF-MCMC simulation (red curve) can obtain similar uncertain range of the QoI as the $f_\mathrm{H}(\mathbf{m})$-MCMC simulation (blue curve). However, using the posterior samples of the $f_\mathrm{L}(\mathbf{m})$-MCMC simulation can result in a wider and more biased prediction of the QoI (black curve).

[Figure 15]

### 3.3. Example 3: Contaminant source identification with 28 unknown parameters

In the first two cases, the numbers of unknown parameters are relatively small. To demonstrate the performance of the AMF-MCMC algorithm in high-dimensional problems, we further test a third numerical case that has 28 unknown parameters. This case is an extension of the second one by considering a more complex contaminant source and heterogeneous conductivity field.

Here the contaminant source is characterized by eight parameters, i.e., the source location $(x_s, y_s)$, and the source strengths $s_i$ during $t = i:i + 1[\text{T}]$, for $i = 1, \ldots, 6$. Here we assume that our prior knowledge about the eight source parameters are rather limited and represented by a multivariate uniform distribution (Table 4).

[Table 4]

The log-transformed conductivity $Y = \ln K$ is assumed to be spatially correlated in the following form:

$$C_Y(x_1, y_1; x_2, y_2) = \sigma_Y^2 \exp\left(-\frac{|x_1 - x_2|}{\lambda_x} - \frac{|y_1 - y_2|}{\lambda_y}\right), \quad (22)$$

where $(x_1, y_1)$ and $(x_2, y_2)$ are two arbitrary locations in the flow domain, $\sigma_Y^2 = 0.4$ is the variance of the $Y$ field, $\lambda_x = 10[\text{L}]$ and $\lambda_y = 5[\text{L}]$ are the correlation lengths in the $x$ and $y$ directions, respectively. To reduce the dimensionality of the $Y$ field, we adopt the Karhunen-Loève (KL) expansion (*Zhang & Lu*, 2004) in this case:

$$Y(\mathbf{x}) \approx \bar{Y}(\mathbf{x}) + \sum_{i=1}^{N_{\text{KL}}} \sqrt{\tau_i} s_i(\mathbf{x}) \xi_i, \quad (23)$$

where $\bar{Y} = 2$ is the mean value of the $Y$ field, $\tau_i$ and $s_i(\mathbf{x})$ are the eigenvalues and eigenfunctions of the correlation function defined in equation (22), $\xi_i \sim \mathcal{N}(0,1)$ are independent Gaussian random variables, for $i = 1, \ldots, N_{\text{KL}}$, respectively. Here $N_{\text{KL}} = 20$ are the number of truncated KL terms, which can preserve about 88.3% of the field variance. Then the unknown model parameters for the heterogeneous $Y$ field are transformed to the 20 KL terms, i.e., $\xi_i$, for $i = 1, \ldots, 20$.

To infer the 28 unknown parameters for the contaminant source and conductivity field, we collect measurements of concentration and hydraulic head at wells denoted by the blue circles in Figure 12. The concentration measurements are collected at $t = 4, 6, 8, 10, 12$[T] and the hydraulic head measurements are collected only once. The measurements are generated from one set of true model parameters and perturbed with measurement errors $\varepsilon_c \sim \mathcal{N}(0, 0.005^2)$ and $\varepsilon_h \sim \mathcal{N}(0, 0.005^2)$. True values of the eight contaminant source parameters are listed in Table 4.

In this case, $f_H(\mathbf{m})$ is the numerical model built with MODFLOW and MT3DMS, while $f_L(\mathbf{m})$ is a data-driven surrogate built with a famous machine learning method, i.e., artificial neural network (ANN), using 2,000 sets of $f_H(\mathbf{m})$ simulation data (80% data are used for training, 10% are used for validation and the rest 10% are used for testing). Here, the Neural Net toolbox in MATLAB R2014a is used. The ANN has 30 hidden layers and is trained with the Levenberg-Marquardt algorithm (*Hagan & Menhaj*, 1994). By comparing $f_L(\mathbf{m})$ and $f_H(\mathbf{m})$ outputs at 200 prior samples, we can obtain $R^2 = 0.960$. Then we apply the $f_H(\mathbf{m})$-MCMC, the $f_L(\mathbf{m})$-MCMC and the AMF-MCMC algorithms to infer the 28 unknown parameters, respectively. In both the $f_H(\mathbf{m})$- and $f_L(\mathbf{m})$-MCMC simulations, the number of parallel chains is set as 20, and the length of each chain is 12,000, which means 240,000 model evaluations in total for each simulation. In the AMF-MCMC simulation, the initial number of $f_H(\mathbf{m})$ evaluations is 100, the initial number of $f_L(\mathbf{m})$ evaluations is 200, the iteration number is 21, and at each iteration 5 sets of $f_H(\mathbf{m})$ data and 5 sets of $f_L(\mathbf{m})$ data are added to refine the GP system locally. In Figure 16, we plot the thinned chains of the $f_H(\mathbf{m})$-MCMC simulation for the eight contaminant source parameters. It is clear that the Markov chains need about 100,000 $f_H(\mathbf{m})$ evaluations to converge to the true parameter values (black crosses). While the results of the $f_L(\mathbf{m})$-MCMC simulation are significantly biased, although we have used 2,000 $f_H(\mathbf{m})$ evaluations to construct the low-fidelity model. With just 205 more $f_H(\mathbf{m})$ and 305 more $f_L(\mathbf{m})$ evaluations, a great improvement of the parameter estimation can be obtained by the AMF-MCMC algorithm.

[Figure 16]

## 4. Conclusions and discussions

In this paper, we propose an efficient method for posterior exploration of hydrologic systems, i.e., the adaptive multi-fidelity MCMC (AMF-MCMC) algorithm. In the AMF-MCMC algorithm, data from both a high-fidelity model and a low-fidelity model are simultaneously fused to build a multi-fidelity GP system, based on which the MCMC simulation can be implemented quickly. As the region of interest is the posterior distribution, we successively add new parameter samples that are close to this region for the high- and low-fidelity model evaluations, which are then used to refine the multi-fidelity GP system locally. Finally, we can obtain an accurate estimation of the posterior distribution with a small number of the high-fidelity model evaluations.

To demonstrate the performance of the AMF-MCMC algorithm, we test three numerical cases in inverse modeling of hydrologic systems. Different types of low-fidelity models are used for illustration. In the first example, we estimate soil hydraulic and thermal parameters with the proposed method in a single ring infiltration experiment. Here the low-fidelity model is built with HYDRUS-2D with a coarser discretization. In the second example, we test a contaminant source identification problem that has multimodal posterior. Here the low-fidelity model is a data-driven surrogate. The third example handles a similar problem as the second one but considers more unknown parameters. In the three examples, the AMF-MCMC algorithm can obtain almost identical results as the $f_\text{H}(\mathbf{m})$-MCMC algorithm but with a very low computational cost. Furthermore, the adaptive multi-fidelity framework is universal in that it can be straightforwardly combined with other inverse or data assimilation methods, e.g., different variants of Kalman filter/smoother (*Chen & Oliver*, 2012; *Emerick & Reynolds*, 2013; *Evensen*, 2007; *Gu & Oliver*, 2007).

When the number of unknown parameters is large (e.g., $N_\text{m} > 100$), many more training data from the high-fidelity model simulations are needed. Then the CPU time of model evaluations will be prohibitive. Moreover, the computational cost in GP system construction will also be huge as it scales cubically with the number of training data. Although we can apply some advanced GP methods, e.g., sparse GP (*Quinonero-Candela & Rasmussen*, 2005; *Snelson & Ghahramani*, 2006), to alleviate the computational cost of GP system construction, the AMF-MCMC algorithm proposed in this paper might still not be a good choice. In this situation, a more computationally appealing strategy is to adopt a method that assumes the posterior distribution to be multi-Gaussian (*Elshall et al.*, 2014; *Rasmussen et al.*, 2015) to roughly estimate the uncertainty. Such assumption is flawed when the posterior is strongly non-Gaussian or even multimodal, then people may take extra techniques like normal-score transform (*Zhou et al.*, 2011) or local updating strategy (*Zhang et al.*, 2018) for a remedy.

# Acknowledgments

Computer codes and data used in this paper are available at:

https://www.researchgate.net/publication/322714837_MATLAB_Codes_of_the_AMF-MCMC_Algorithm

This work is supported by the National Natural Science Foundation of China (grants 41771254 and 41571215).

The authors would like to thank the editor and anonymous reviewers for their constructive comments and suggestions, which significantly improve the quality of this work. The authors would also like to thank Andreas Klimke from University of Stuttgart for providing the sparse grid interpolation toolbox, Maziar Raissi from Brown University for providing the MATLAB codes of multi-fidelity GP, Jasper Vrugt from University of California, Irvine for providing the MATLAB code of DREAM$_{(ZS)}$, respectively.

# Tables

Table 1 Prior ranges and true values of unknown model parameters in the first example

| Parameters | Prior ranges | True values |
|---|---|---|
| $\alpha [\text{cm}^{-1}]$ | [0.0190  0.0930] | 0.0387 |
| $n [-]$ | [1.360  2.370] | 2.210 |
| $K_s [\text{cm h}^{-1}]$ | [4.828  11.404] | 6.759 |
| $b_1 [\text{kg cm h}^{-3} \text{K}^{-1}]$ | $[2.179 \times 10^{12}\ \ 4.857 \times 10^{12}]$ | $2.948 \times 10^{12}$ |
| $b_2 [\text{kg cm h}^{-3} \text{K}^{-1}]$ | $[4.778 \times 10^{11}\ \ 2.426 \times 10^{12}]$ | $2.118 \times 10^{12}$ |
| $b_3 [\text{kg cm h}^{-3} \text{K}^{-1}]$ | $[2.174 \times 10^{12}\ \ 5.184 \times 10^{12}]$ | $2.972 \times 10^{12}$ |

Table 2 Values of known model parameters in the first example

| Parameters | Values | Parameters | Values |
|---|---|---|---|
| $\theta_r [\text{cm}^3 \text{cm}^{-3}]$ | 0.041 | $l [-]$ | 0.500 |
| $\theta_s [\text{cm}^3 \text{cm}^{-3}]$ | 0.430 | $C_w [\text{J cm}^{-3} \text{K}^{-1}]$ | 4.180 |
| $\theta_n [\text{cm}^3 \text{cm}^{-3}]$ | 0.600 | $C_n [\text{J cm}^{-3} \text{K}^{-1}]$ | 1.920 |
| $\theta_o [\text{cm}^3 \text{cm}^{-3}]$ | 0.001 | $C_o [\text{J cm}^{-3} \text{K}^{-1}]$ | 2.510 |
| $\lambda_x [\text{cm}]$ | 0.200 | $\lambda_y [\text{cm}]$ | 2.000 |

Table 3 Prior ranges and true values of unknown model parameters in the second example

| Parameters | $x_s [\text{L}]$ | $y_s [\text{L}]$ | $S_s [\text{MT}^{-1}]$ | $t_{on} [\text{T}]$ | $t_{off} [\text{T}]$ |
|---|---|---|---|---|---|
| Ranges | [3  5] | [3  7] | [10  13] | [3  5] | [9  11] |
| True values | 3.854 | 5.999 | 11.044 | 4.897 | 9.075 |

Table 4 Prior ranges and true values of contaminant source parameters in the third example

| Parameter | Range | True value |
|---|---|---|
| $x_s[L]$ | [3 5] | 4.033 |
| $y_s[L]$ | [4 6] | 5.405 |
| $s_1[MT^{-1}]$ | [0 8] | 1.229 |
| $s_2[MT^{-1}]$ | [0 8] | 7.628 |
| $s_3[MT^{-1}]$ | [0 8] | 4.327 |
| $s_4[MT^{-1}]$ | [0 8] | 5.438 |
| $s_5[MT^{-1}]$ | [0 8] | 0.293 |
| $s_6[MT^{-1}]$ | [0 8] | 6.474 |

# Figures

Figure 1. Comparison between the single-fidelity Gaussian process (SF-GP) simulation (magenta dashed curve) and the multi-fidelity Gaussian process (MF-GP) simulation (black dashed curve). Here the $f_H(m)$ and $f_L(m)$ simulations are represented by the blue curve and red curve, the training data for SF-GP and MF-GP are represented by the magenta squares and black circles, respectively.

Figure 2. Flow domain for the first example.

Figure 3. Comparison of simulated (a) head and (b) temperature outputs between $f_H(\mathbf{m})$ and $f_L(\mathbf{m})$.

Figure 4. Trace plots of the $\hat{R}$-statistic of the six model parameters in (a) $f_H(\mathbf{m})$-MCMC simulation and (b) $f_L(\mathbf{m})$-MCMC simulation. The threshold of 1.2 for convergence diagnosis is represented by the red dashed lines.

Figure 5. Trace plots of model parameters obtained by the $f_H(\mathbf{m})$-MCMC algorithm (red dots) and the $f_L(\mathbf{m})$-MCMC algorithm (blue dots) in the first example. The true values are represented by the black crosses.

Figure 6. Evolution of (a-b) variance of the multi-fidelity GP output and (c-d) RMSE between the multi-fidelity GP output and the high-fidelity output averaged over 400 posterior samples.

Figure 7. $N_H$ initial parameter samples (red dots) and $I_{\max}$ successively added parameter samples (blue dots) for $f_H(\mathbf{m})$ evaluations in the AMF-MCMC algorithm. The true values are represented by the black crosses.

Figure 8. Marginal PPDFs obtained by the $f_H(\mathbf{m})$-MCMC algorithm (red curves) and the AMF-MCMC algorithm (blue dashed curves). The true values are represented by the black vertical lines.

Figure 9. Finally obtained (a-b) variance of the multi-fidelity GP output and (c-d) RMSE between the multi-fidelity GP output and the high-fidelity output when adding

different numbers of new training data sets at each iteration.

Figure 10. Trace plots of model parameters evaluated by $f_H(\mathbf{m})$ in the AMF-MCMC algorithm (blue dots) and the AGP-MCMC algorithm (red dots). The true values are represented by the black crosses.

Figure 11. Marginal PPDFs obtained by the $f_H(\mathbf{m})$-MCMC algorithm (red curves), the AMF-MCMC algorithm (blue dashed curves), and the AGP-MCMC algorithm (magenta dash-dotted curves), respectively. The true values are represented by the black vertical lines.

Figure 12. Flow domain for the second and third examples.

Figure 13. Trace plots of model parameters obtained by the $f_H(\mathbf{m})$-MCMC algorithm (red dots) and the $f_L(\mathbf{m})$-MCMC algorithm (blue dots) in the second example. The true values are represented by the black crosses.

Figure 14. Bivariate scatter plots of posterior parameter samples obtained by the $f_H(\mathbf{m})$-MCMC algorithm (red dots) and the AMF-MCMC algorithm (blue dots) in the second example. The true values are represented by the black crosses.

Figure 15. Predictive uncertainty of the QoI based on estimation results of the $f_H(\mathbf{m})$-MCMC simulation (blue curve), the $f_L(\mathbf{m})$-MCMC simulation (black curve) and the AMF-MCMC simulation (red curve), respectively. The true value of the QoI is represented by the black cross.

Figure 16. Trace plots of the $f_H(\mathbf{m})$-MCMC simulation (blue dots), posterior mean estimates obtained by the $f_L(\mathbf{m})$-MCMC simulation (red circles) and the AMF-MCMC simulation (red squares) of the eight contaminant source parameters.

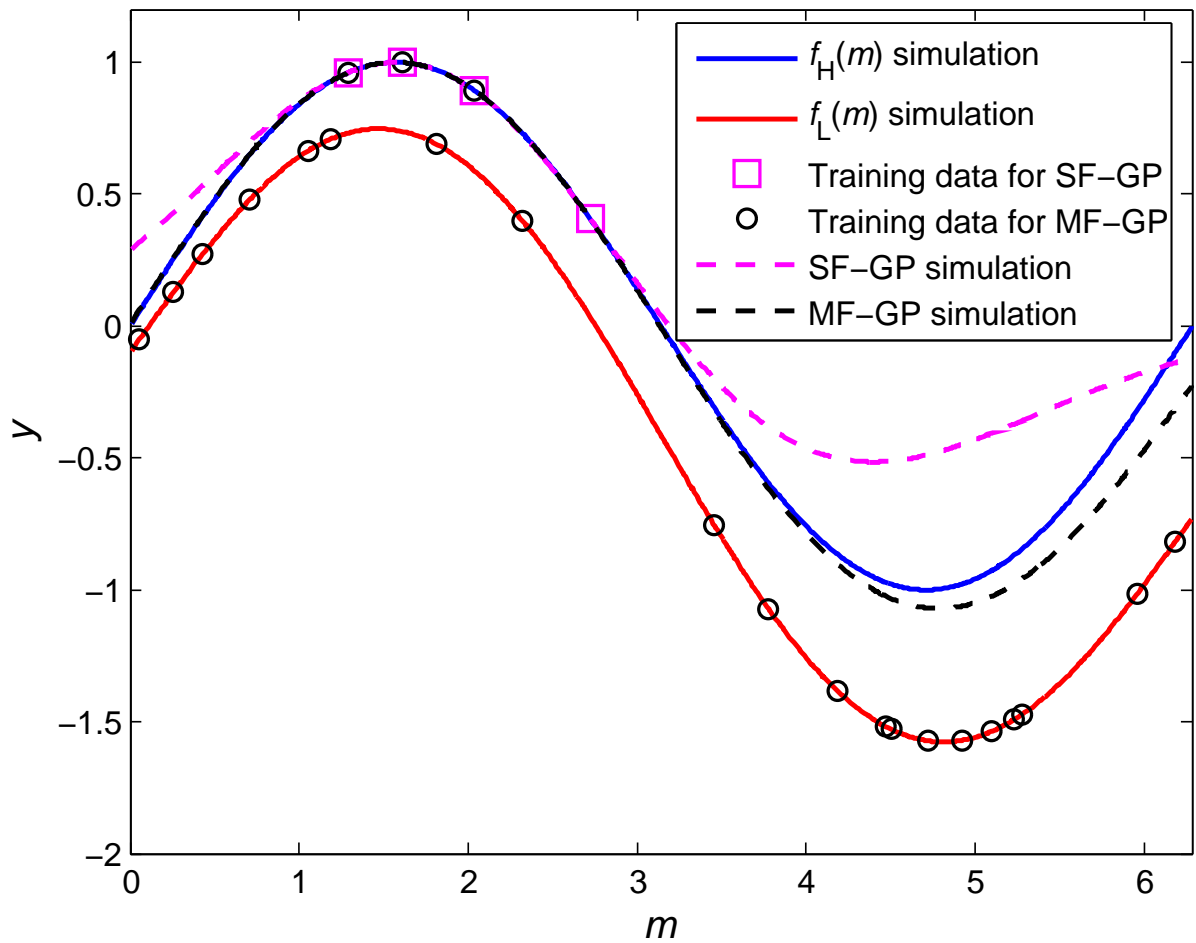

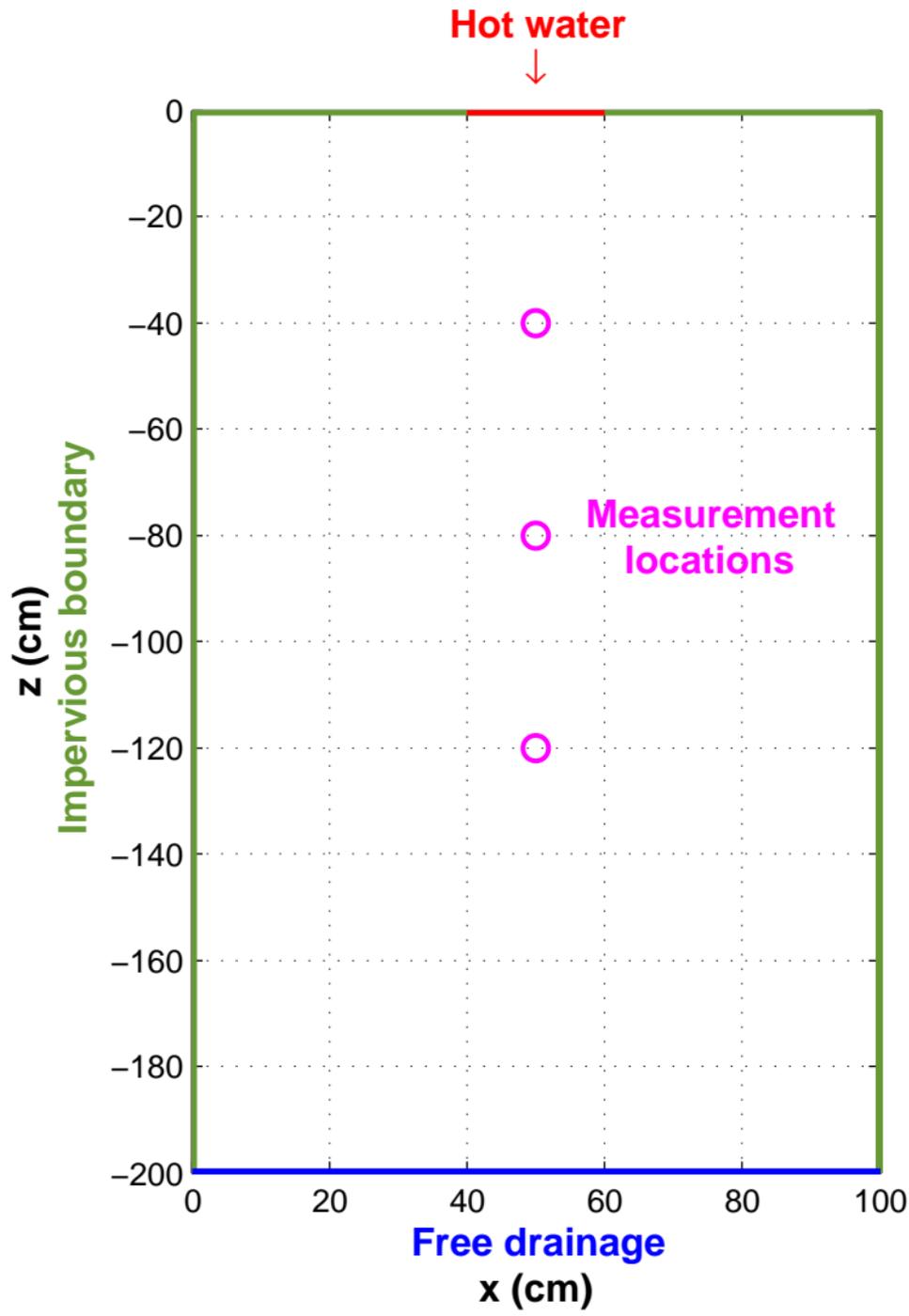

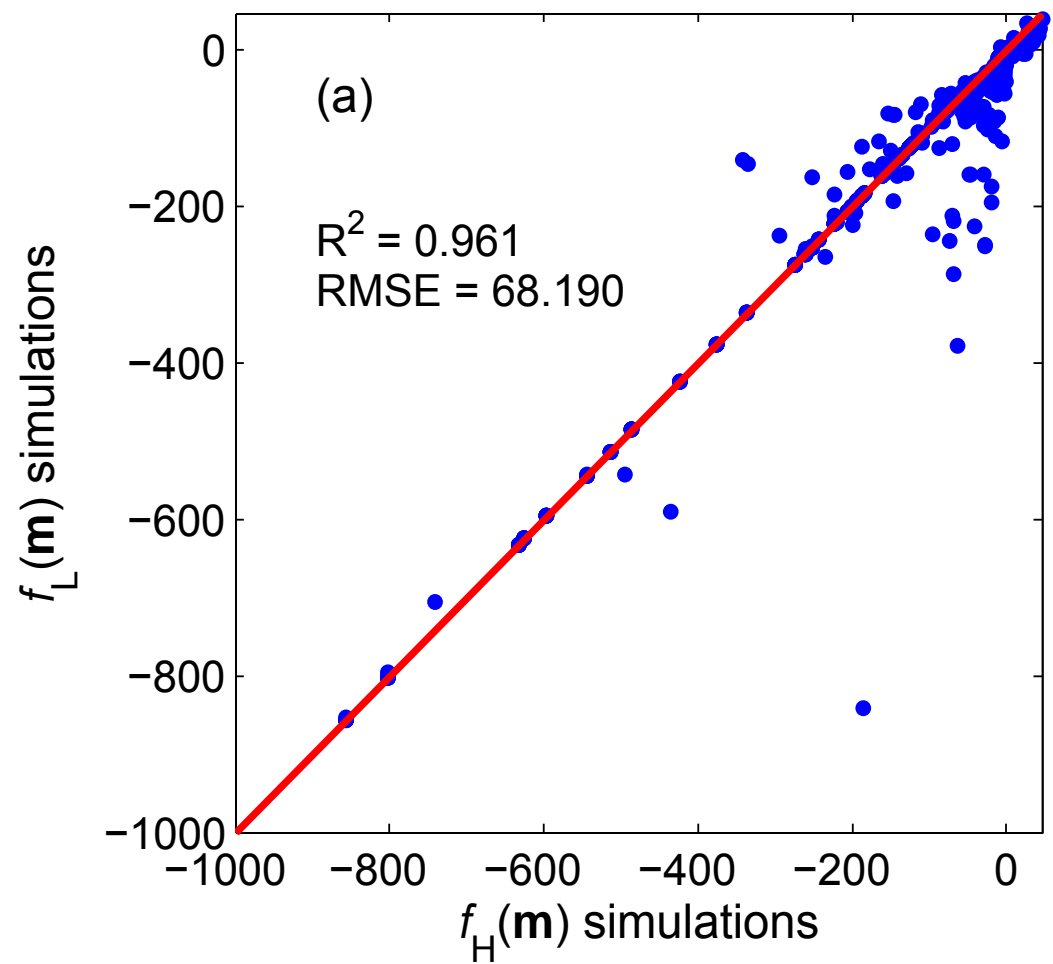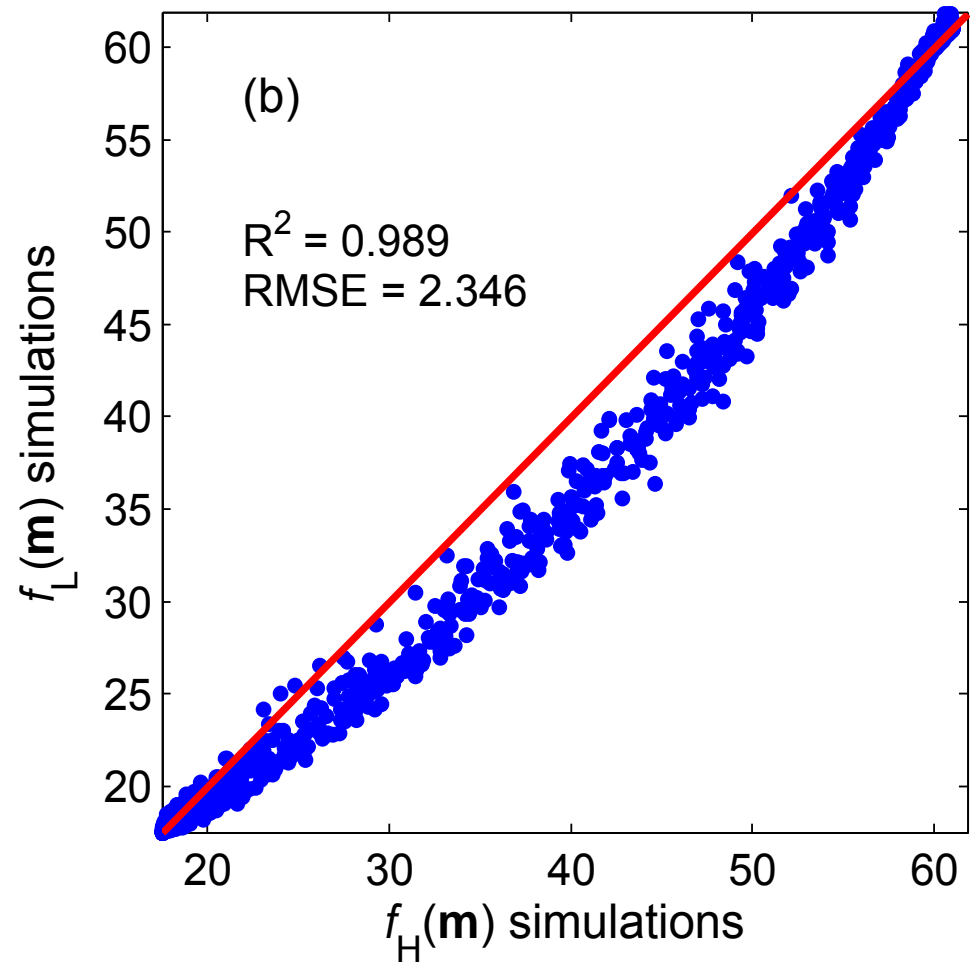

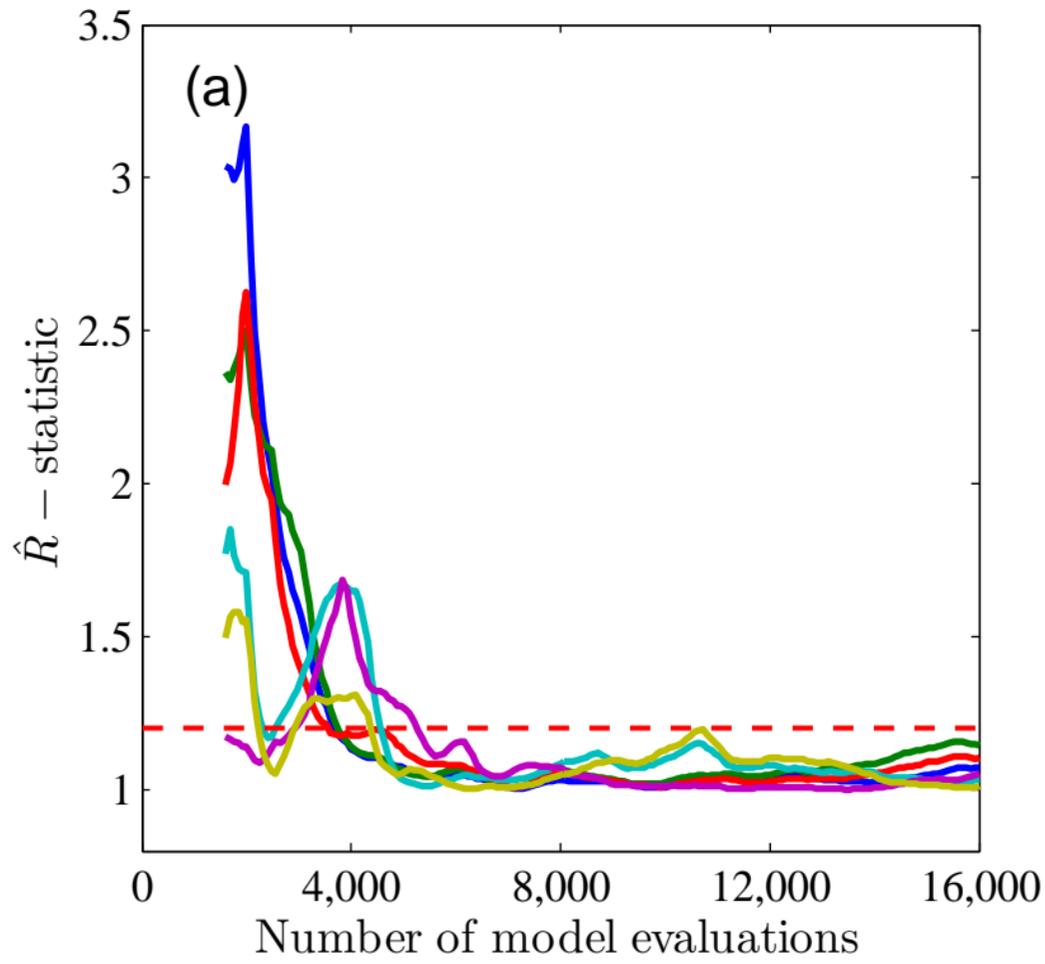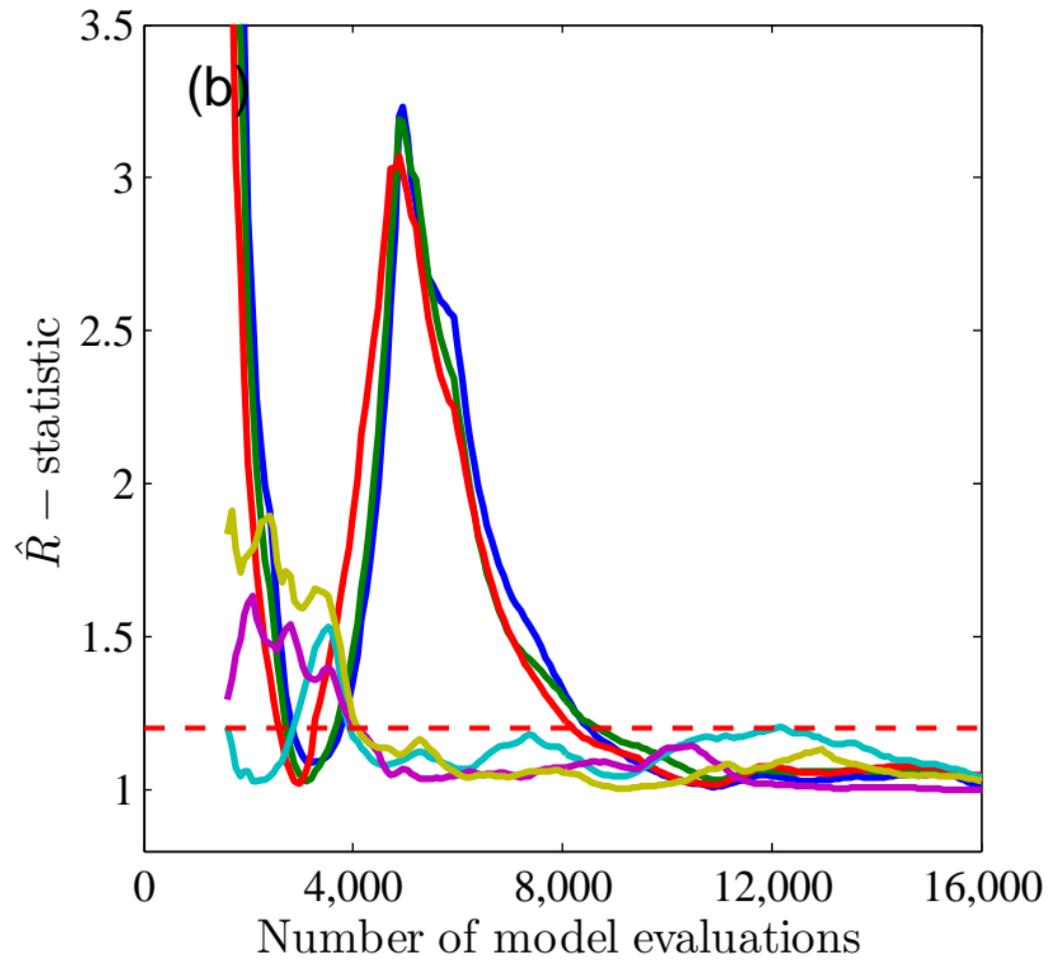

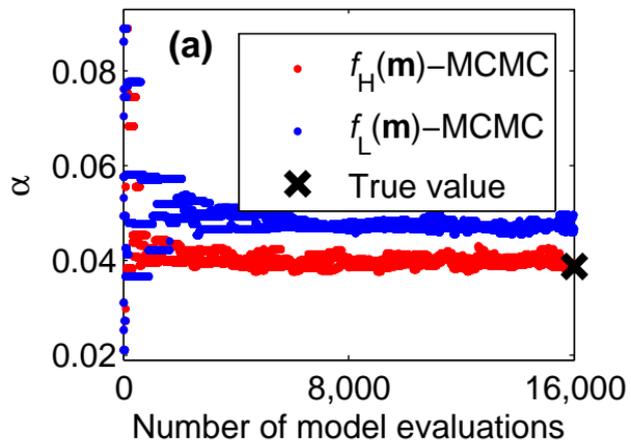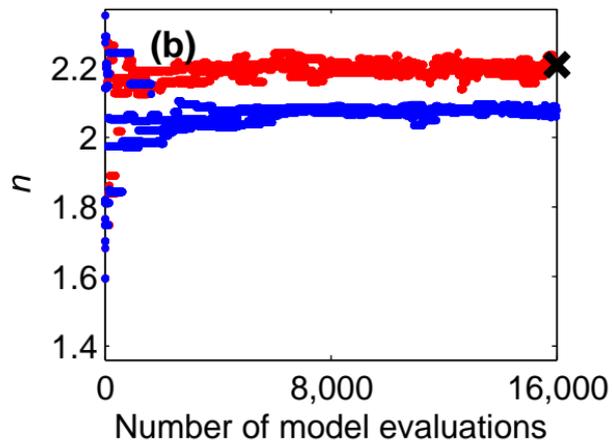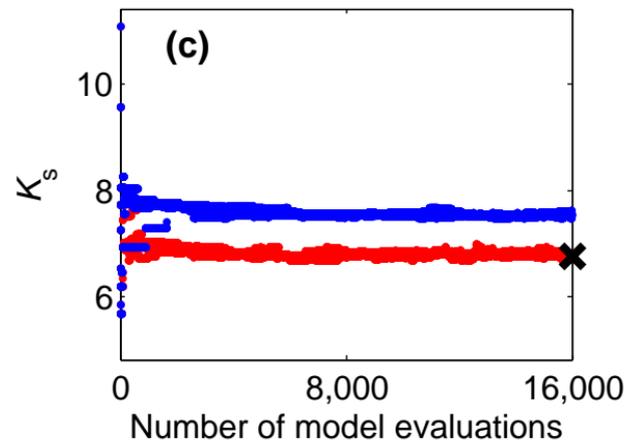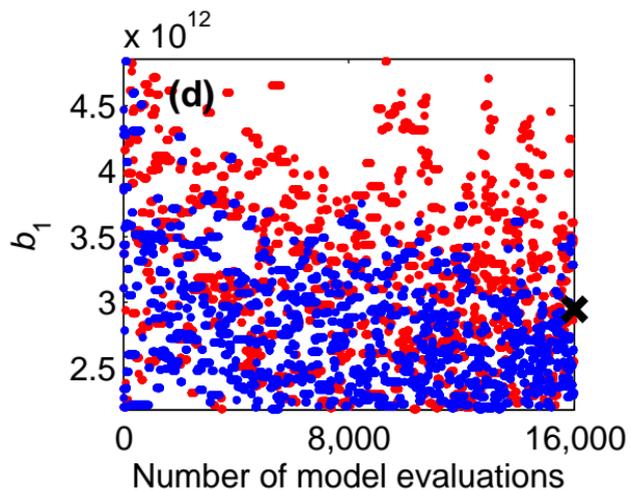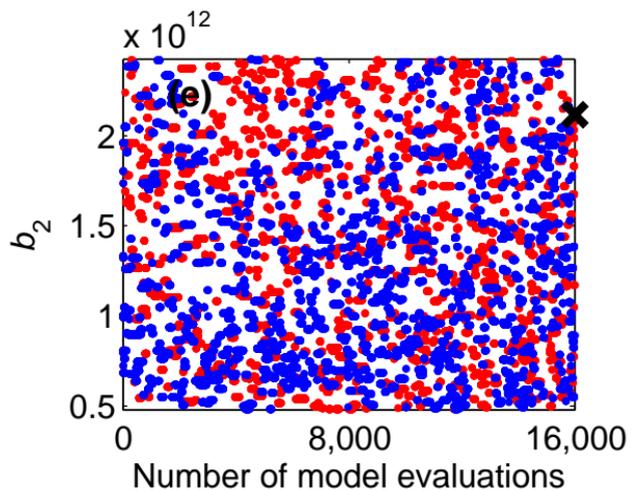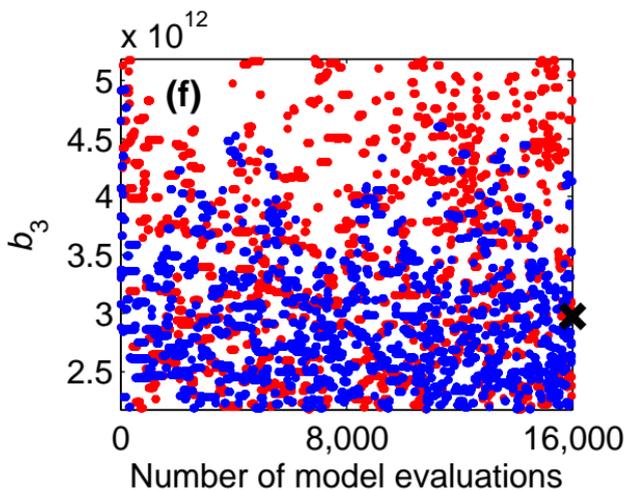

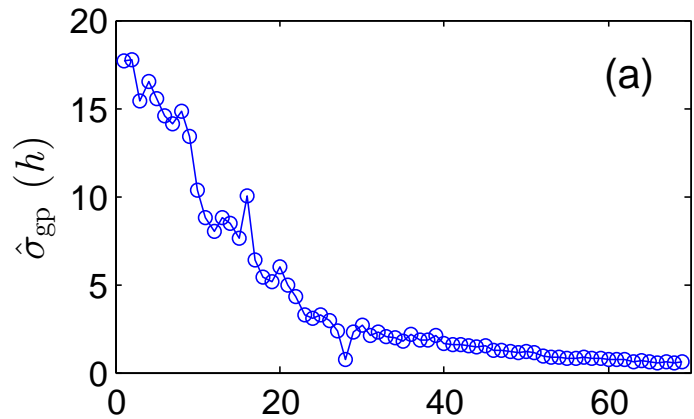
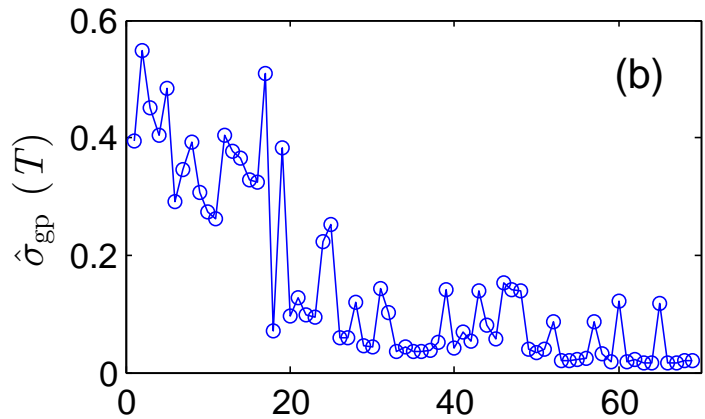
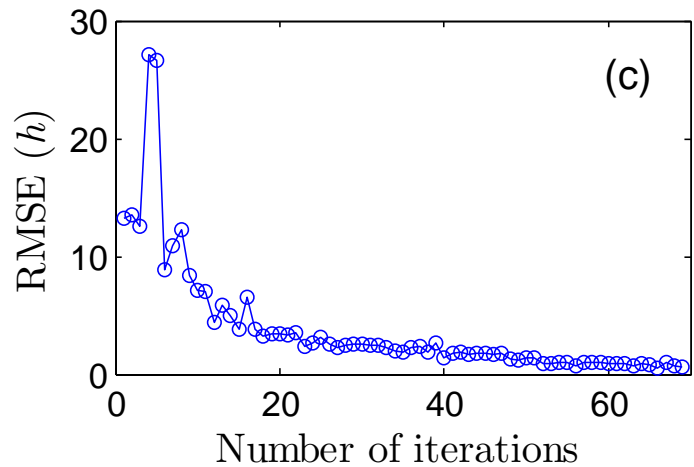
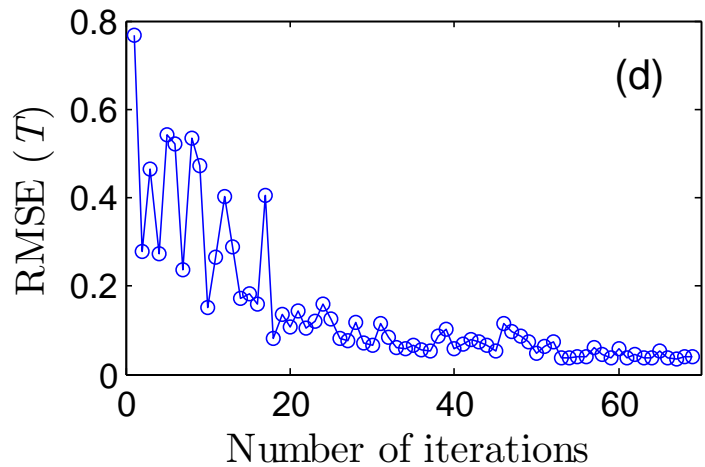

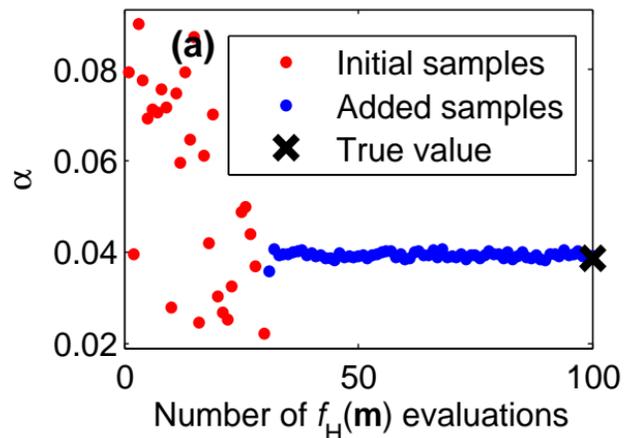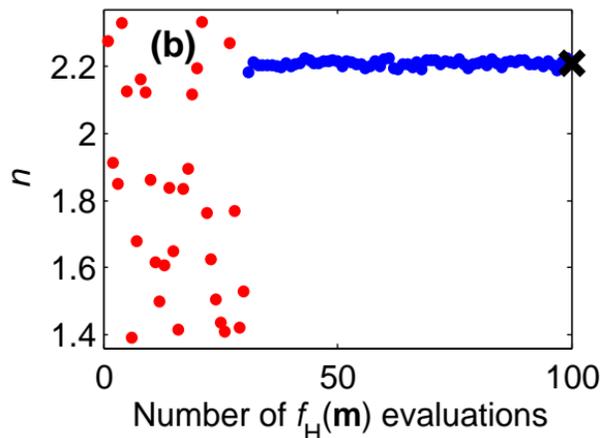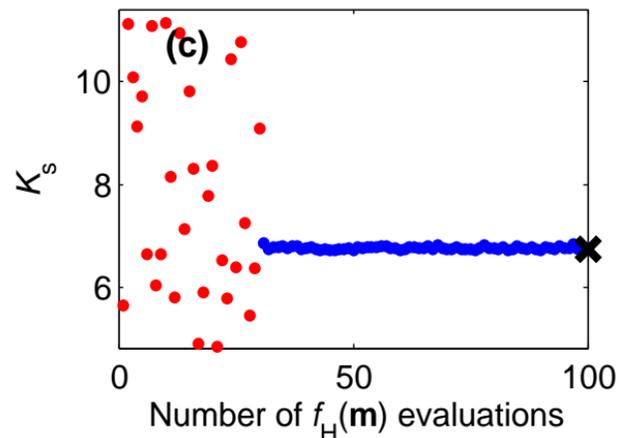
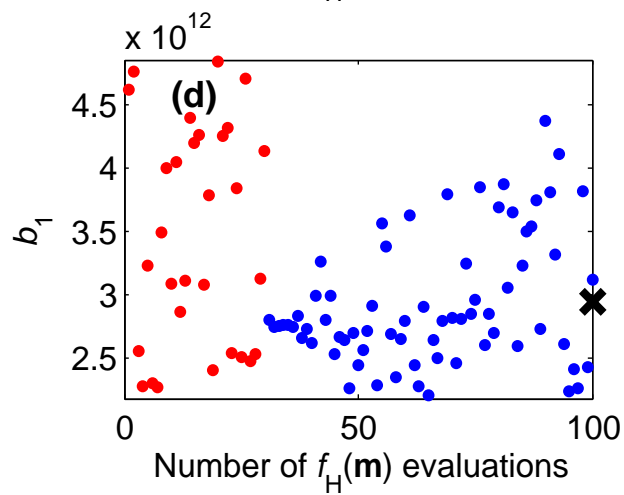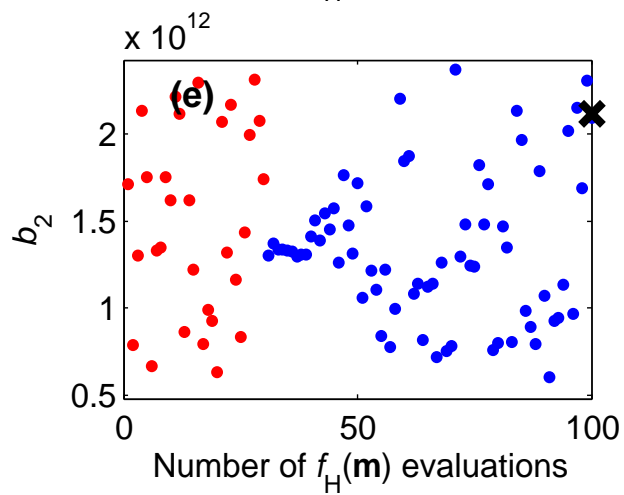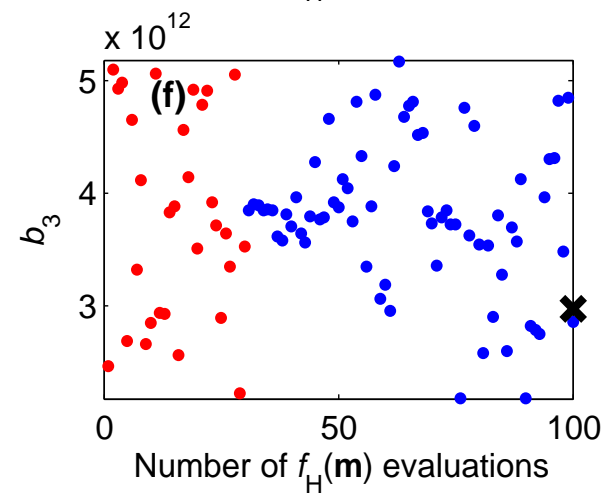

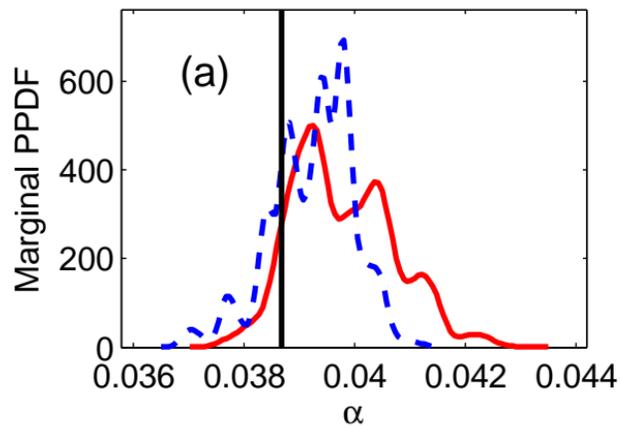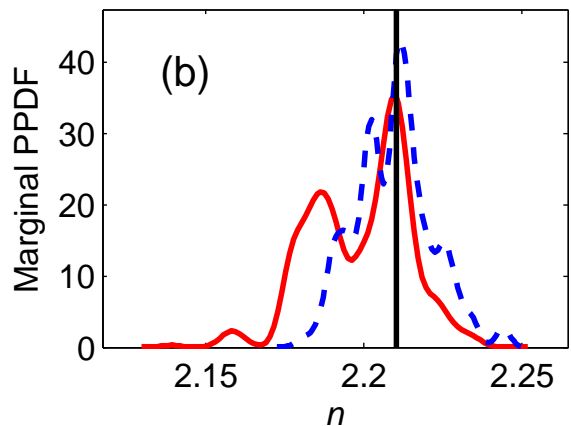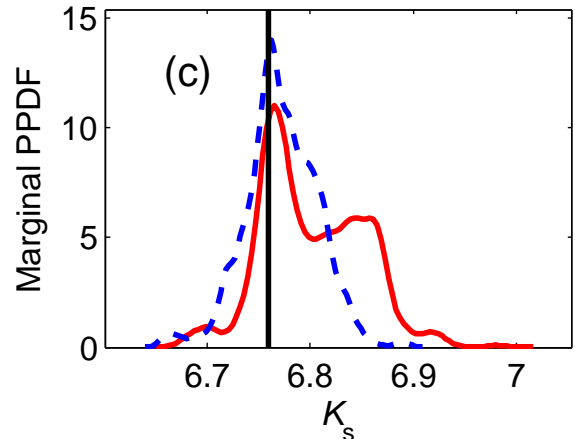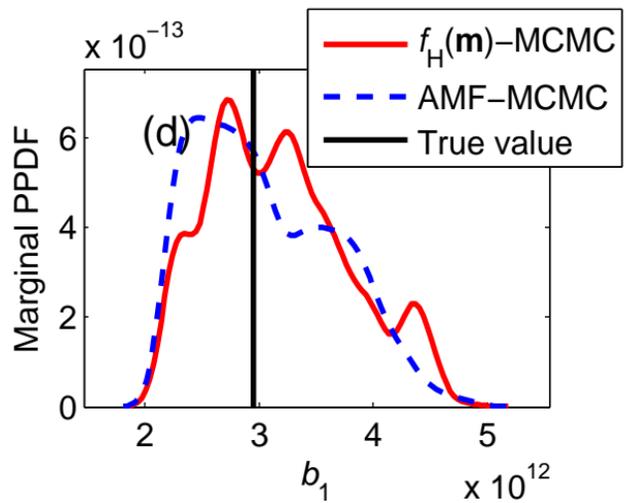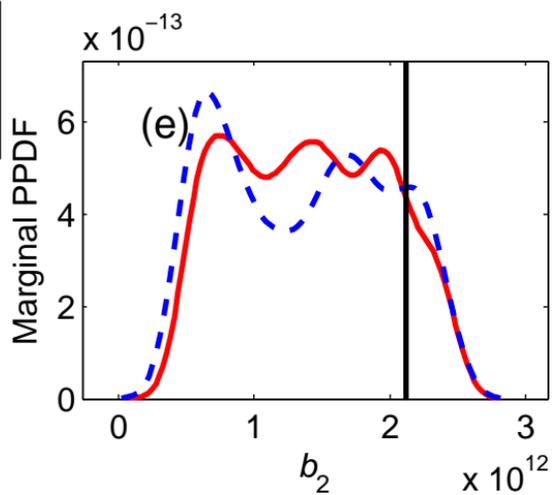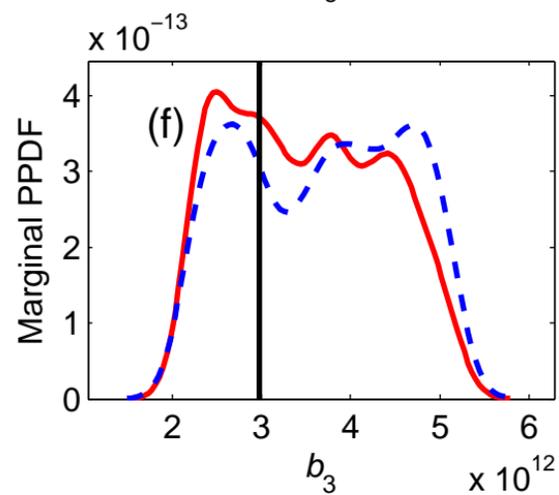

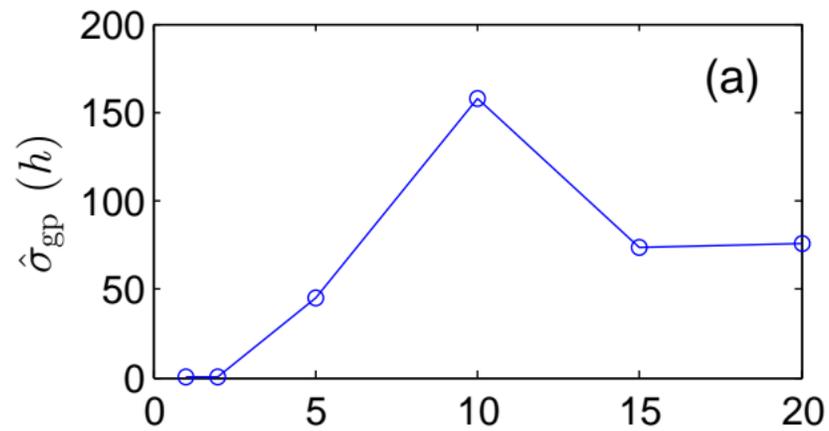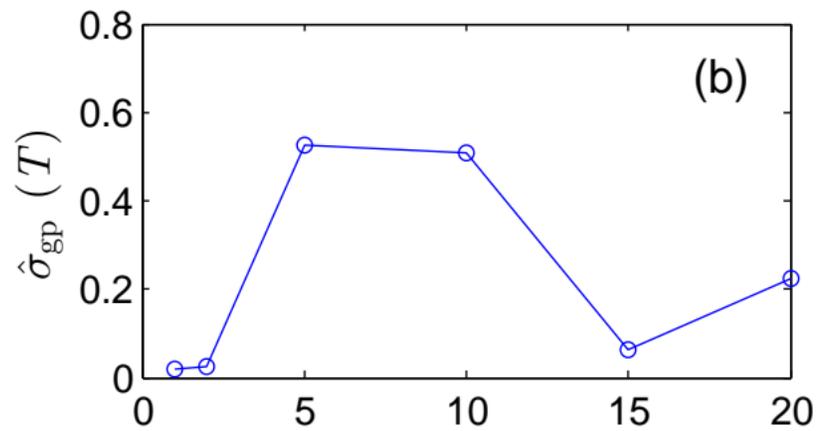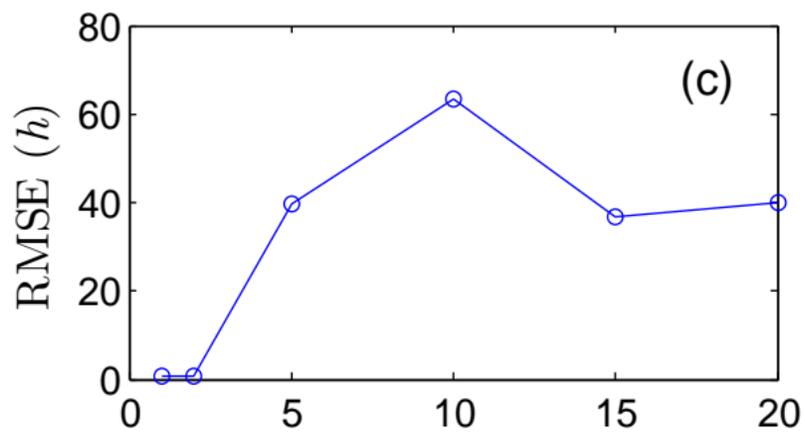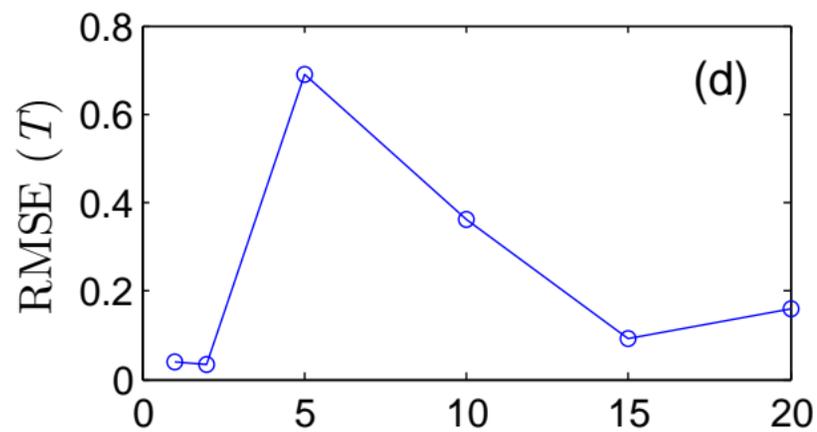

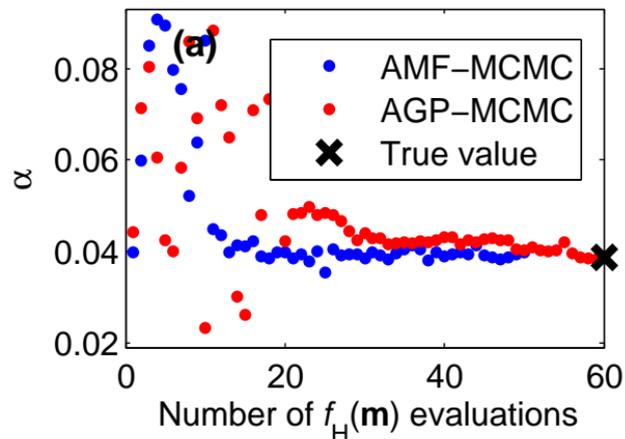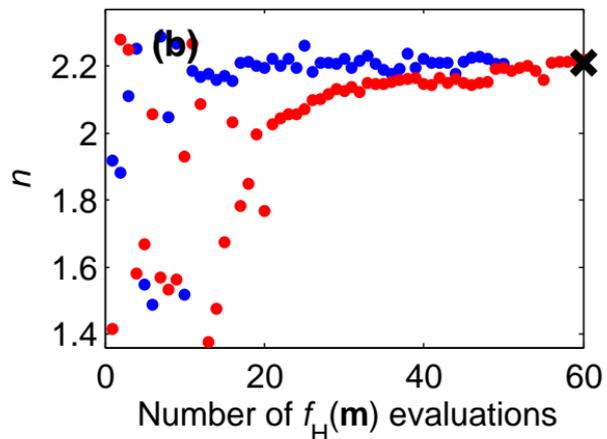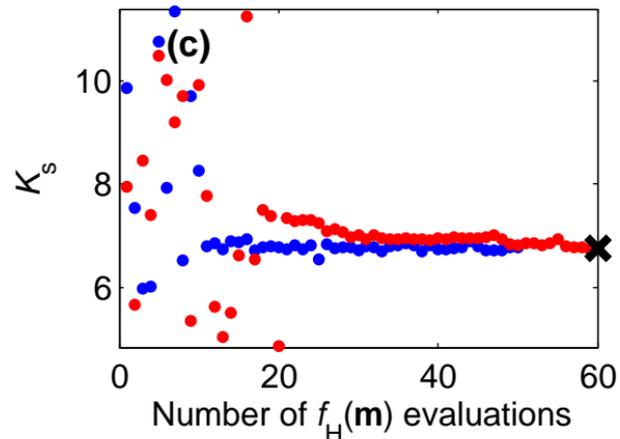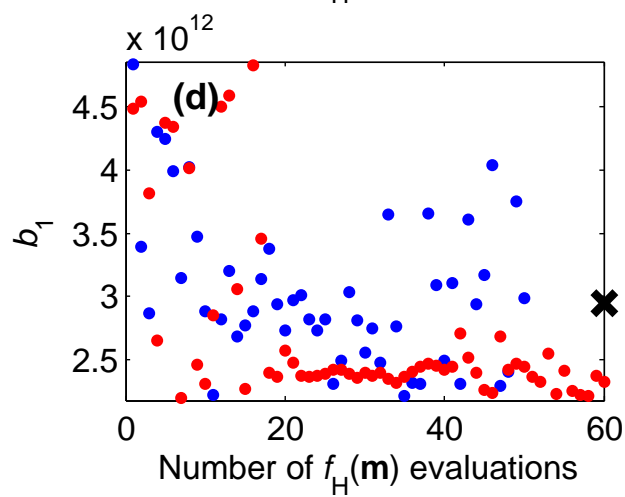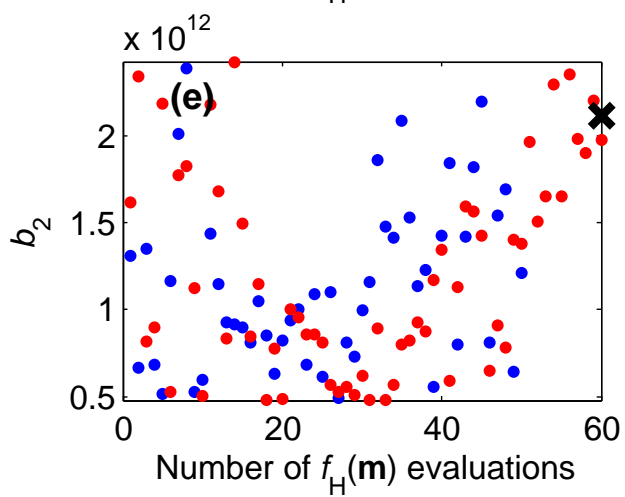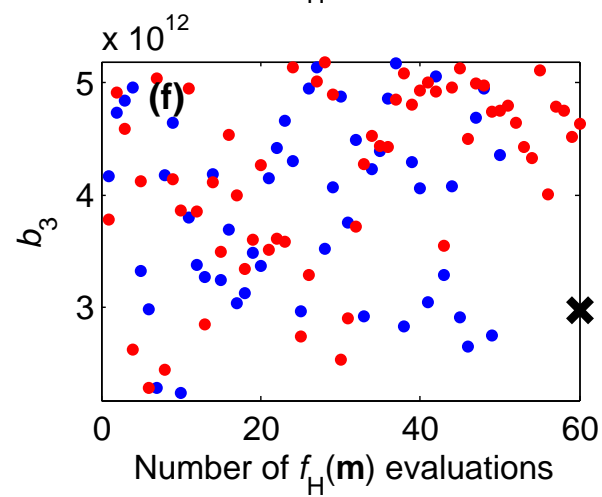

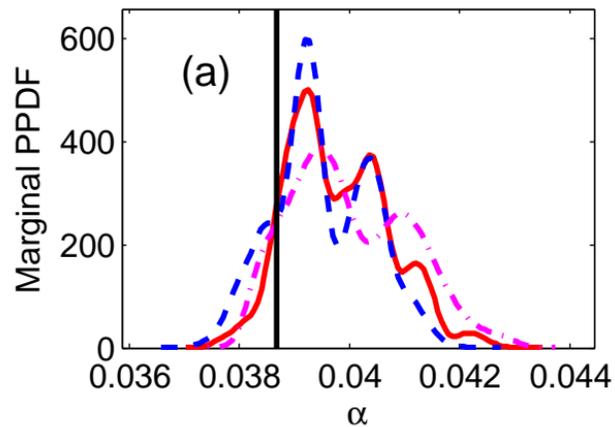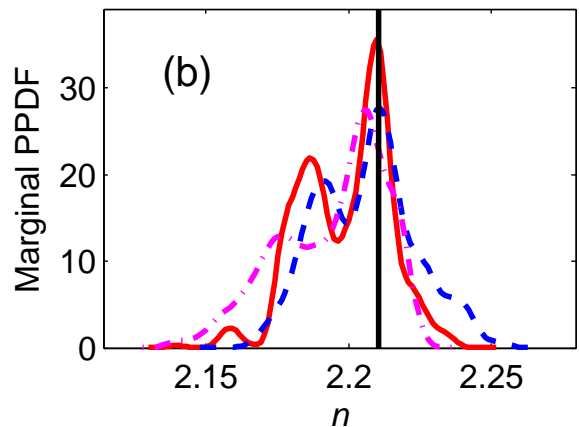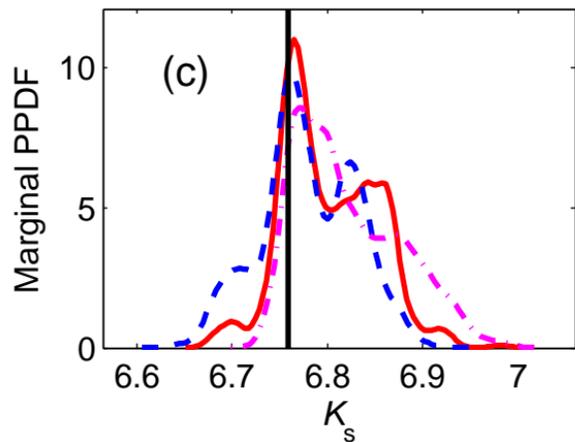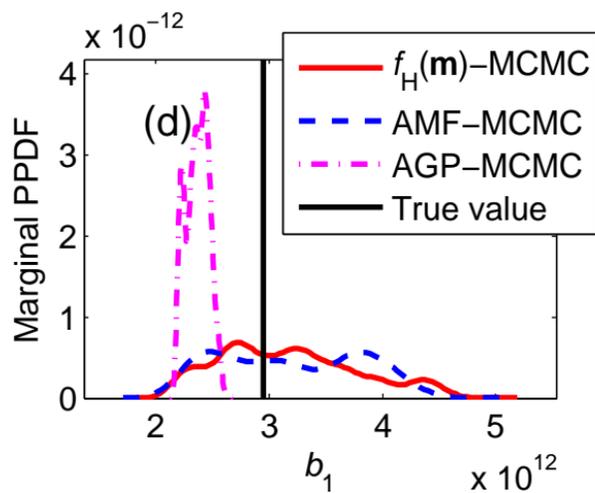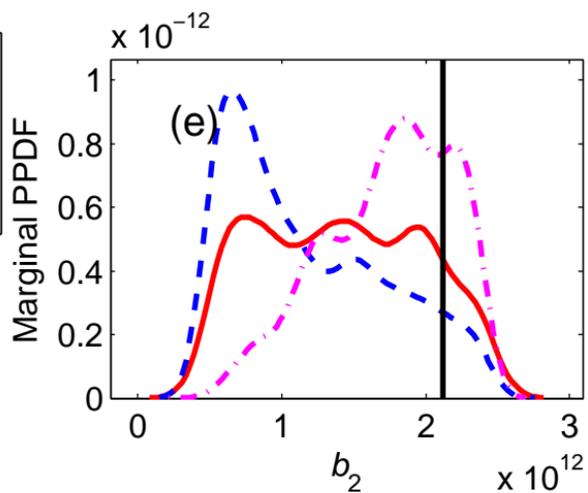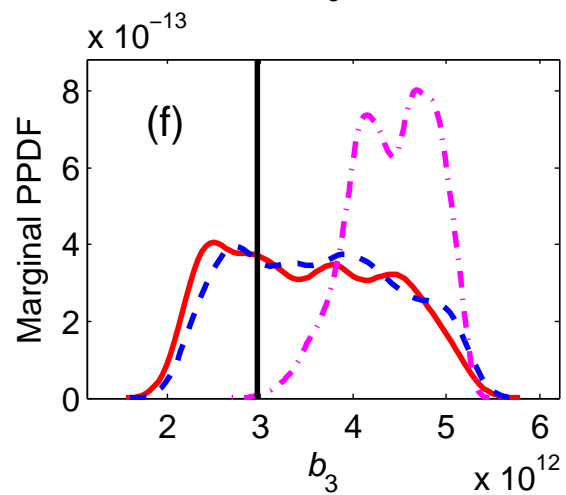

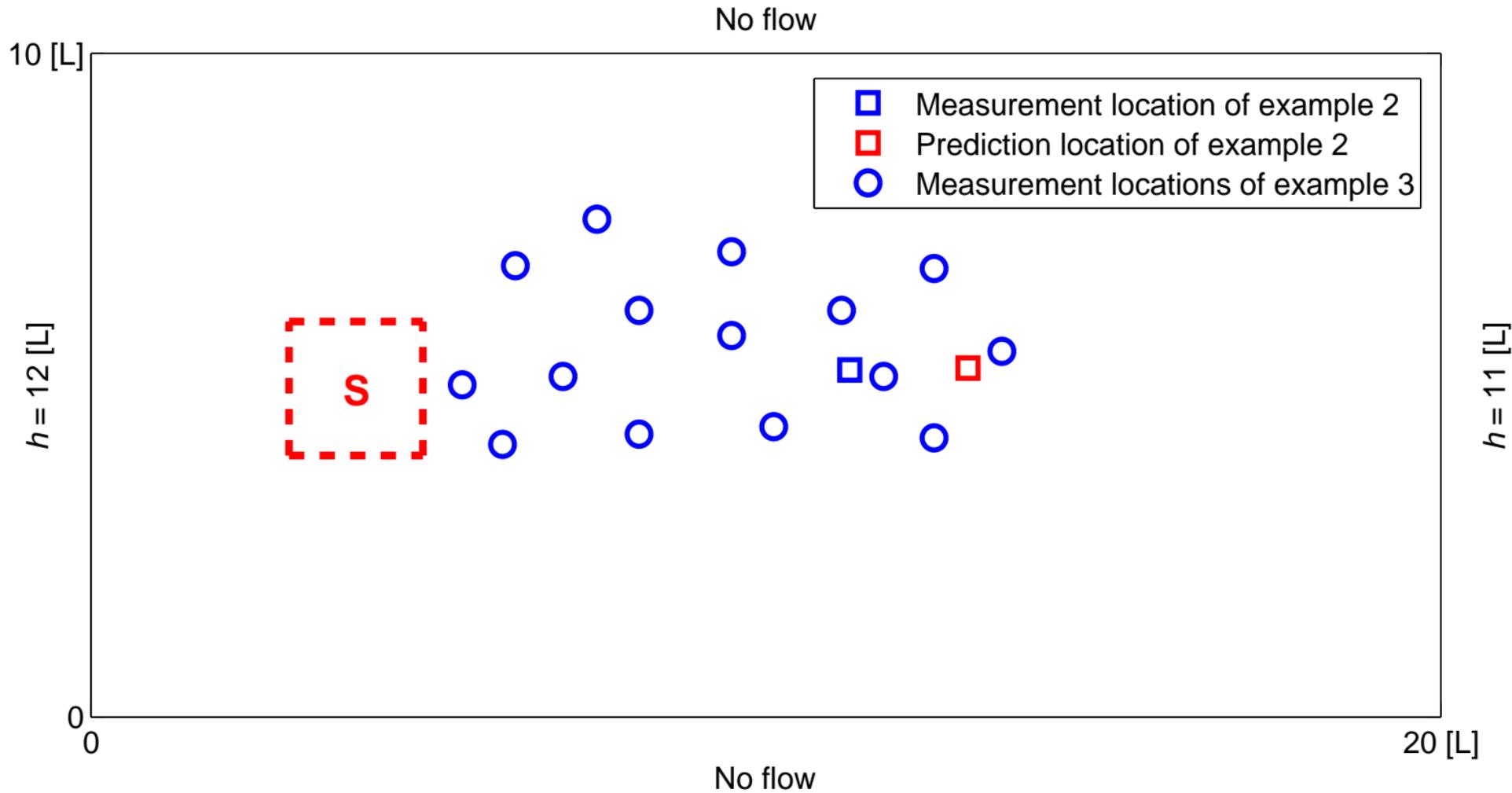

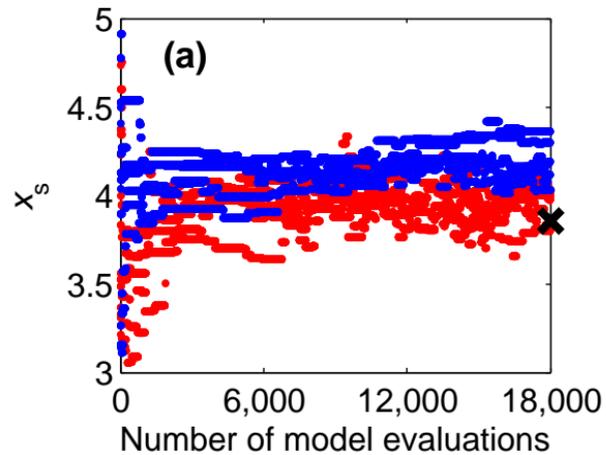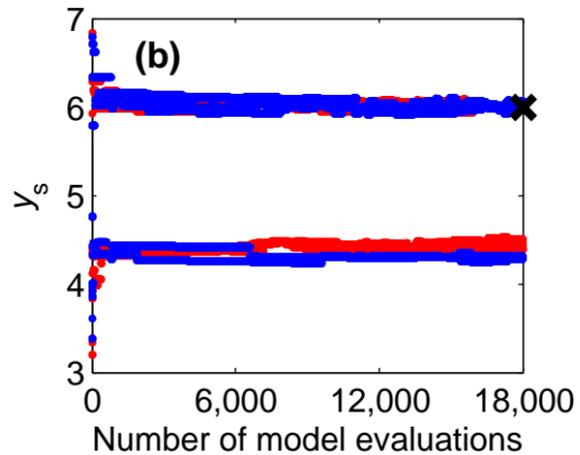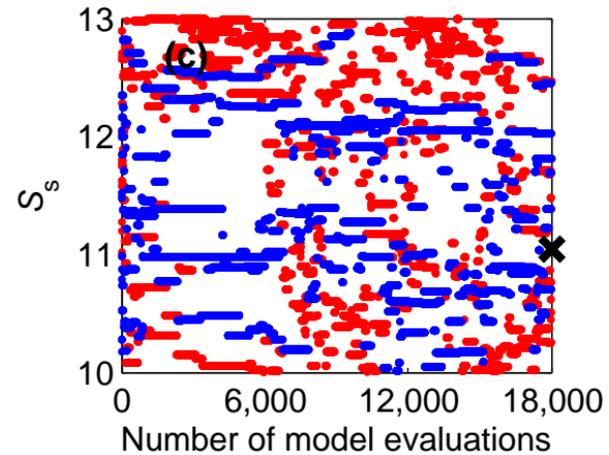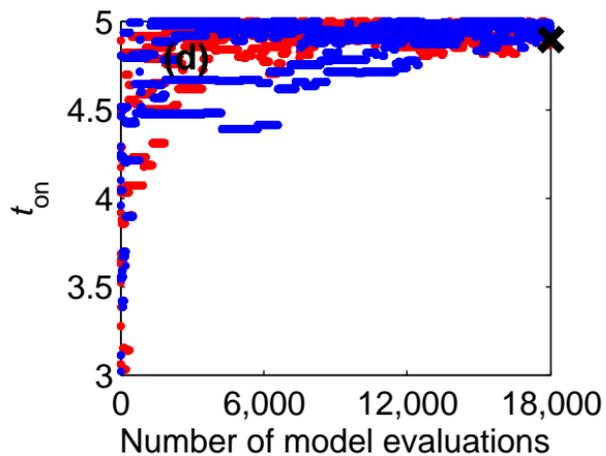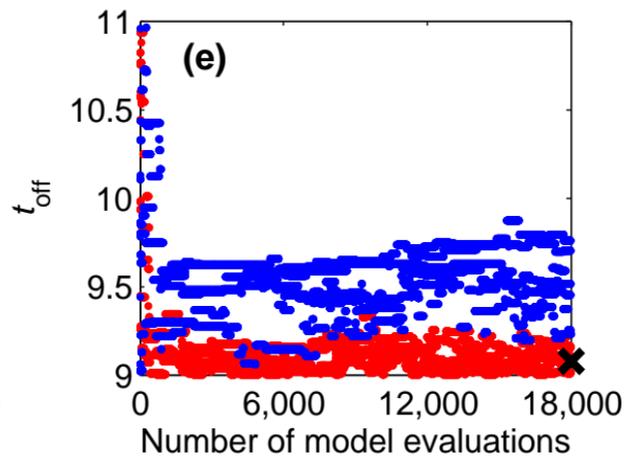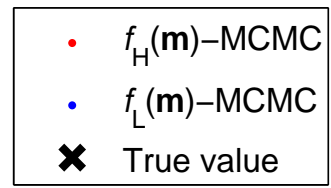

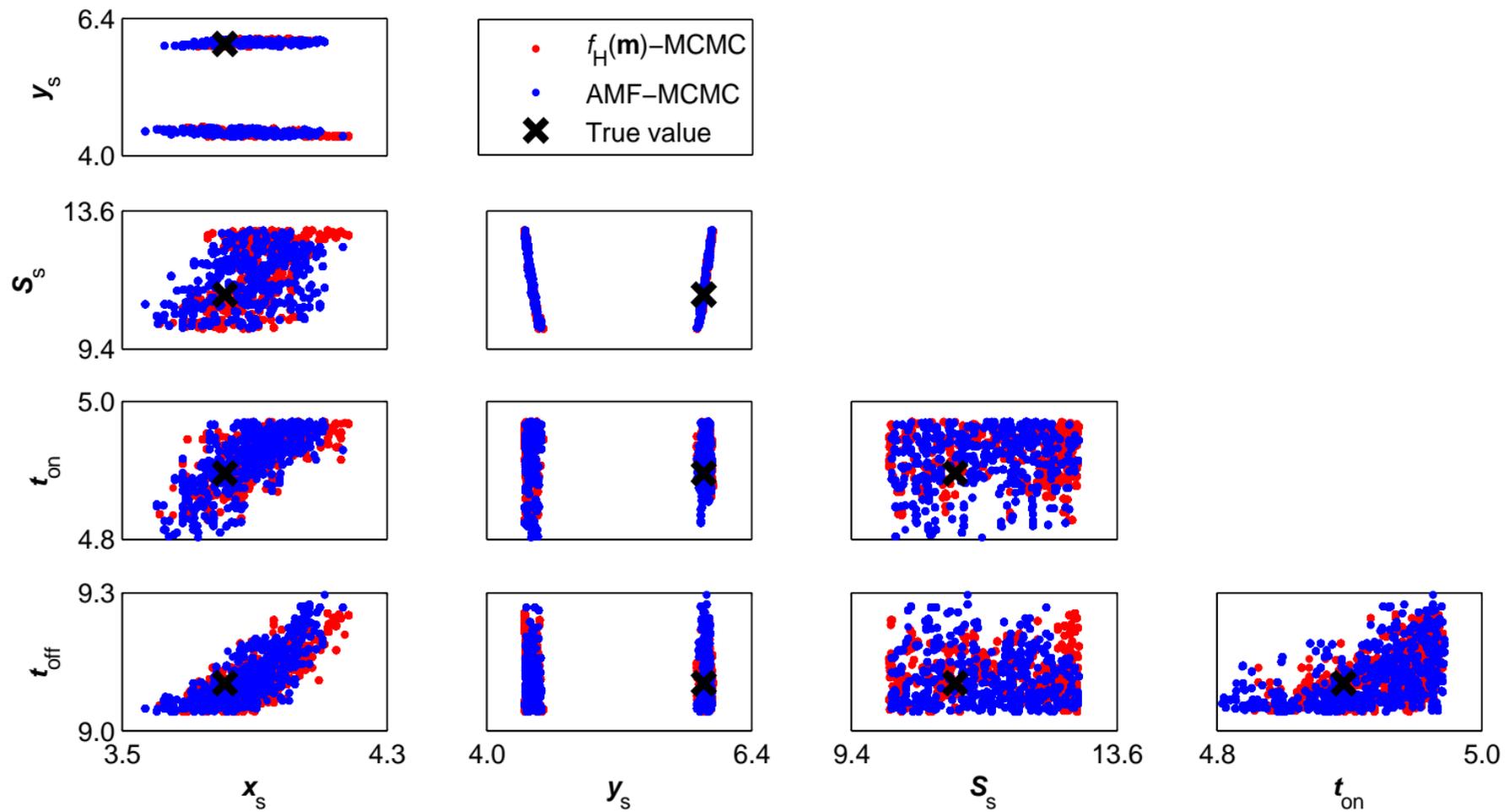

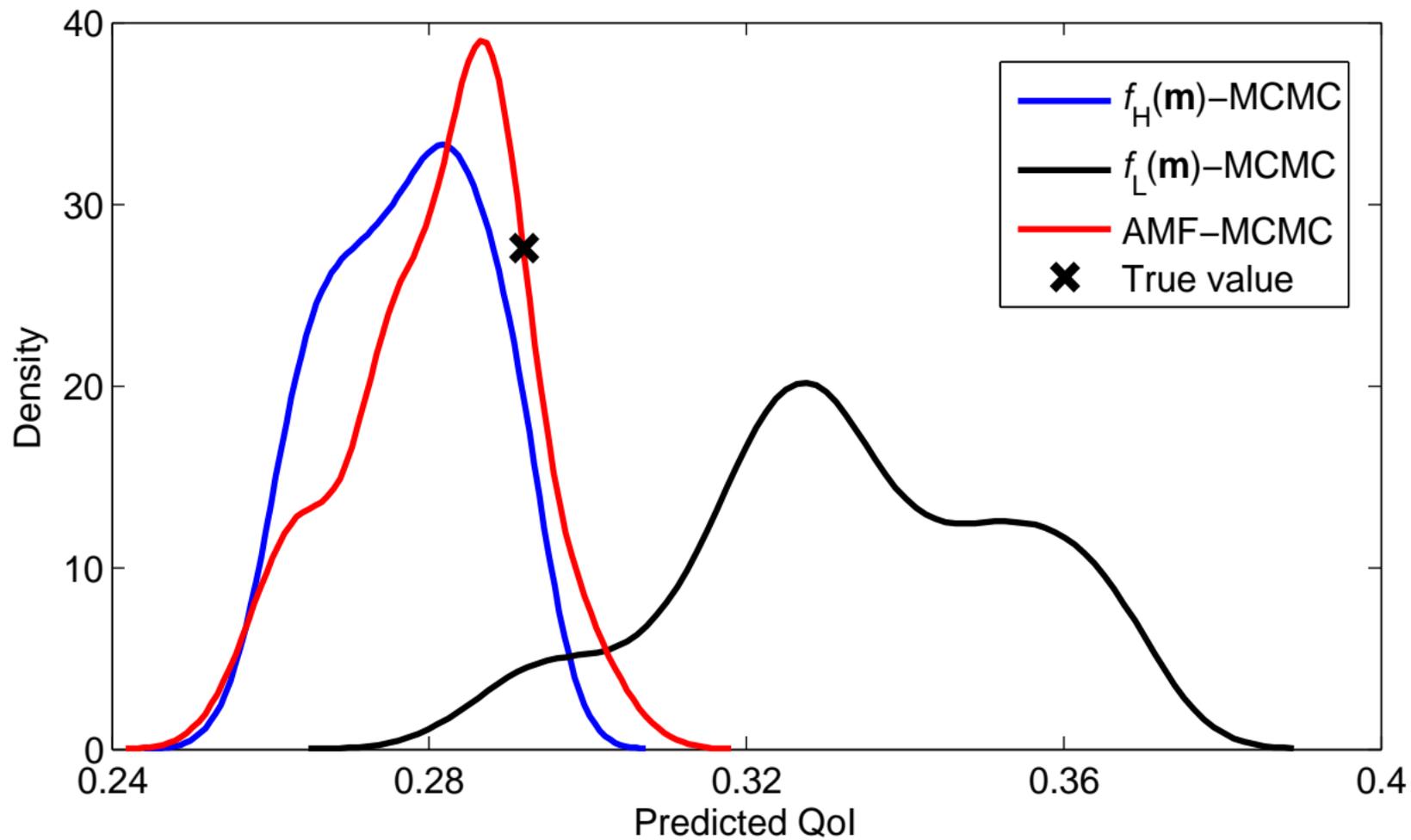

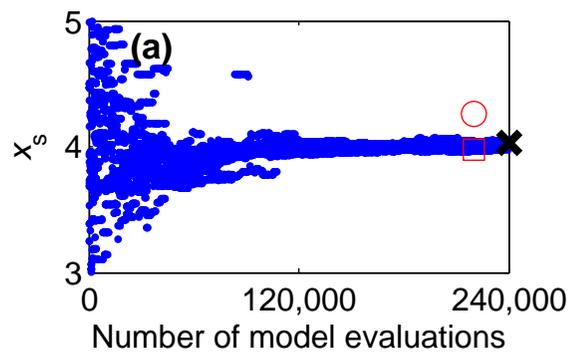
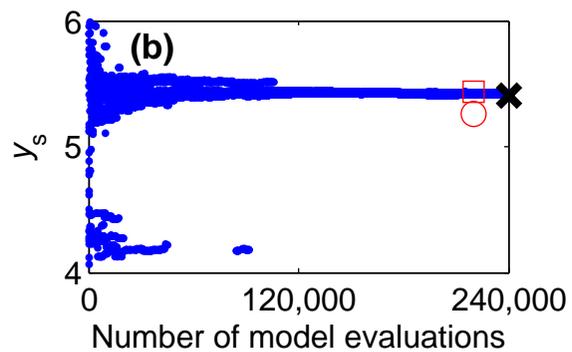
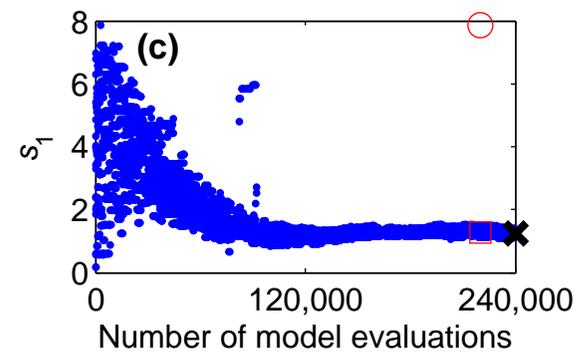
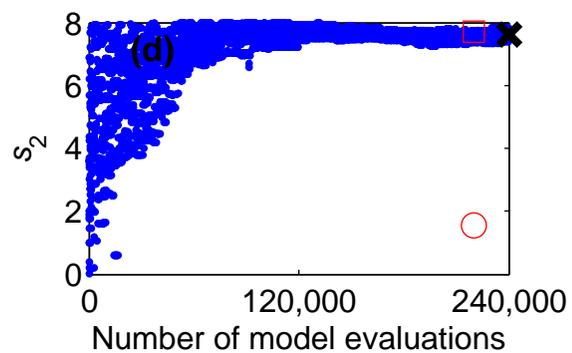
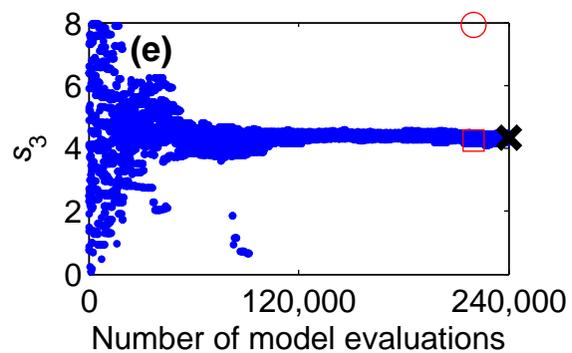
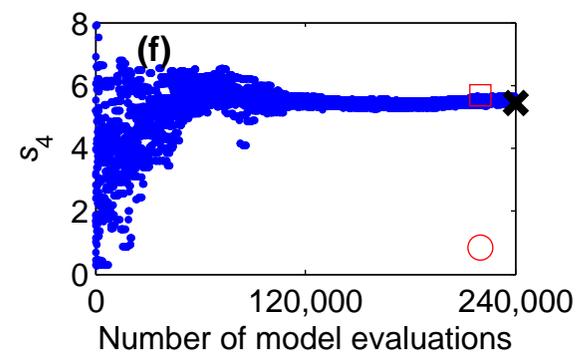
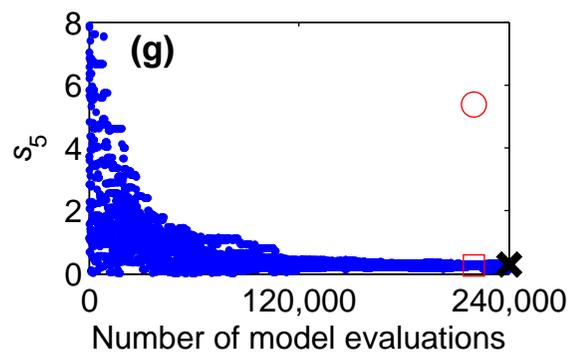
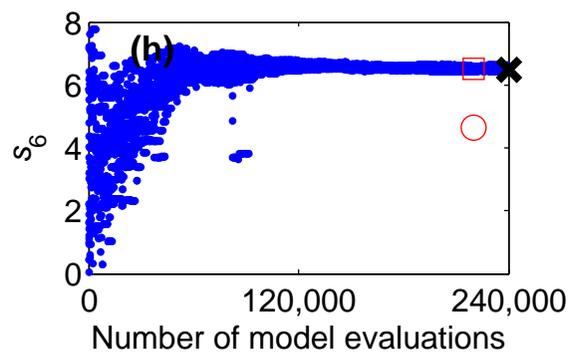
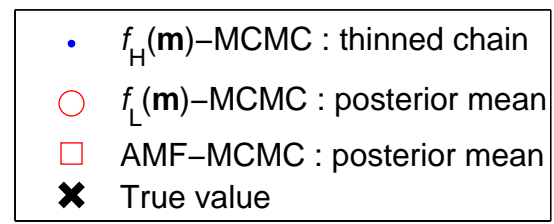